\documentclass[11pt]{article}

\usepackage{amsmath,amssymb}
\usepackage{amsthm}
\usepackage{graphicx}
\usepackage{hyperref}
\usepackage[nottoc,numbib]{tocbibind}
\usepackage{geometry}
\newtheorem{theorem}{Theorem}[section]
\newtheorem{lemma}[theorem]{Lemma}

\newtheorem{corollary}[theorem]{Corollary}
\newtheorem{definition}[theorem]{Definition}
\newtheorem{remark}[theorem]{Remark}
\newtheorem{example}[theorem]{Example}

\DeclareMathOperator\var{var}
\DeclareMathOperator\W{\mathcal{W}}

\DeclareMathOperator\argmin{argmin}

\DeclareMathOperator\subg{subg}

\DeclareMathOperator\diam{diam}

\DeclareMathOperator\G{\mathcal{G}}

\begin{document}

\title{Complexity Regularization and Local Metric Entropy}
\author{Fabi\'an Latorre}
\date{2014}
\maketitle 

\newpage
\tableofcontents 

\newpage

\begin{section}{Introduction}\label{Introduction}
\vspace{10 mm}

Machine Learning, a blend  of statistics and computer science, has been an active field of research in recent years. Applications include characters and face recognition, spam filtering and medical diagnosis. In the case of supervised learning one has at hand random sample $(X_i,Y_i)$ from an unknown distribution of a random pair $X \in S$ and $Y \in \mathcal{R}$. Given a class of functions $\mathcal{G}:S \rightarrow \mathbb{R}$ and a risk functional $L: \mathcal{G} \rightarrow \mathbb{R}$, which depends on the distribution of $(X,Y)$, we wish to find $g \in \mathcal{G}$ with minimal risk $L(g)$. Usual risk functions are the absolute expected deviation and the expected square loss.\\

In practice, the task of finding such a function is rather impossible, as the underlying distribution is not known, and one resorts to finding a function that minimizes the empirical risk defined by the random sample, and then finding a confidence interval for the real loss of the function. In the case where $Y$ is a $\{ 0,1 \}$ random variable the problem is known as \emph{classifier selection} and the consistency of such methods have been shown to depend upon a theorem of uniform convergence of relative frequencies of events to their true probabilities, over a class of sets. During the 1970's, Vladimir Vapnik and Alexey Chervonenkis, proved that such a theorem holds if and only if a certain combinatorial quality of the class $\mathcal{G}$, called the VC dimension, is finite. Confidence intervals for the real risk of a classifier, criteria for model selection, among other results can be defined in terms of this quantity, and numerous computations and bounds of such numbers for different classes are readily available in the literature.\\

However, in the case where $Y$ takes values on the interval $[0,1]$, which we refer to as the problem of \emph{Regression Estimation}, consistency of methods based on empirical risk minimization requires uniform convergence of means to expectations, and necessary and sufficient conditions for such convergence become rather difficult to define based on $\mathcal{G}$. Sufficient conditions have been shown to depend on the topology of the class when equipped with the $L_1(P)$ norms, where $P$ is a probability measure on the underlying space. Such conditions measure in some way, the complexity of the class. \\

In the context of Structural Risk Minimization, one is presented a sequence of classes $\{\mathcal{G}_j\}$ from which, given a random sample $(X_i,Y_i)$ one wants to choose a strongly consistent estimator. For certain types of classes of functions, we present a criterion to choose an estimator, based on the minimization of the sum of empirical error and a complexity penalty $r(n,j)$ over each class $\G_j$.\\

We present also several other results together with important results found on current literature on the subject, in an attempt to present and unify the theory in the context of regression estimation. In particular we present a generalization of the consistency of Structural Risk Minimization in the context of regression estimation which are all new results found on Chapter 5.

\end{section}
\newpage

\begin{section}{Classifier Selection}\label{Classifier_Selection}
\vspace{5 mm}

In this section we present the main foundation of the subject. Most of the results are due to Vapnik and Chervonenkis, and most of the proofs can be found, for example, in the book on the subject by Luc Devroye, L\'aszl\'o Gyorfi and G\'abor Lugosi \cite{devroye96}. Throughout we use $S$ to denote a subset of $\mathbb{R}^d$ and we refer to a function $g: S \rightarrow \{0,1\}$ as a classifier .\\

\begin{subsection}{Uniform Convergence of Relative Frequencies to Probabilities}\label{Frequencies}

\begin{definition}  $X_1, \ldots, X_n$, be an i.i.d. random sample from a distribution $P$ on a set $S$. The associated empirical measure $P_n$ is the random measure 
\[
P_n(A) := \dfrac{1}{n} \sum_{i=1}^n 1_{ A}(X_i)
\]
\\
For a function $g: S \rightarrow \mathbb{R}$ we define its empirical integral as
\[
P_nf:=\dfrac{1}{n} \sum_{i=1}^n f(X_i)
\]
As usual denote $P(A)$ the probability of the set $A$ and $P(g)$ the expected value of $g$.
\end{definition}

\begin{definition}\label{Glivenko}Let $P$ be a probability measure on a set $S$. We say that a class of functions $\mathcal{F}$ with domain $S$, is a P-Glivenko-Cantelli Class if 
$$
P \left \{ \sup_{g\in \mathcal{G}} \left |P_n (g) - P(g)\right | > \epsilon \right \} \rightarrow 0 \, \, \mbox{a.s.} \, \, \mbox{ as $n\rightarrow \infty $} 
$$

A class of sets $\mathcal{A}$ is called a \emph{P-Glivenko-Cantelli Class} if its class of characteristic functions $ \{ 1_A : A \in \mathcal{A}\} $ is P-Glivenko-Cantelli. \\
\end{definition}

The supremum on the definition is a minor nuisance as it might not be measurable. Some measurability conditions on the class must be assumed for most of the results presented to hold. In the case of most practical examples such conditions are easily met. All results presented will assume such conditions hold. Refer to $\cite{pollard84}$, Appendix C.\\

\begin{definition}\label{measurability} Let $T$ be a separable metric space and let $\mathcal{B}(T)$ be its Borel $\sigma$-algebra. Let $(S, \mathcal{C})$ be a measurable space and $\mathcal{F}= \{f(\cdot,t) : t \in T \}$ be a class of functions on $S$ indexed by $T$ such that:
\begin{enumerate}
\item The function $f(\cdot, \cdot)$ is $\mathcal{C} \otimes \mathcal{B}(T)$-measurable as a function from $S \otimes T$ into the real line.

\item T is an analytic subset of a compact metric space $\overline{T}$ (from which it inherits its metric and Borel $\sigma$-algebra).
\end{enumerate}
\end{definition}

To prove that a class of functions $\G$ is a $P$-Glivenko-Cantelli class one makes heavy use of the following corollary to the well known Borel-Cantelli Lemma:

\begin{theorem}{(Borel-Cantelli)}\label{BC} For a sequence of random variables $\{ X_n : n=1,2, \ldots \}$ and a random variable $X$ 
\[
\sum_{n=1}^\infty P \left \{  |X_n-X | > \epsilon \right \} < \infty \, \, \, \mbox{ then } \, \, \, X_n \rightarrow X \, \, \, a.s.
\]

\end{theorem}

If a class of functions $\mathcal{G}$ is P-Glivenko-Cantelli we also say that a uniform law of large numbers holds over $\mathcal{G}$. A class of sets $\mathcal{A}$ is known to be P-Glivenko-Cantelli for any probability measure if and only if certain combinatorial quantity that depends on $\mathcal{A}$ is finite.

\begin{definition} Let $\mathcal{A}$ be a class of subsets of $\mathcal{R}^d$ and $E\subset \mathcal{R}^d$ be subset. Let $E \cap \mathcal{A} := \{ E \cap A \,: \, A \in \mathcal{A} \}$ and 
$$
\Delta_\mathcal{A}(n):= \sup_{ \{E : |E|=n \} } |E \cap \mathcal{A}|
$$
When $\Delta_\mathcal{A}(n)$ is bounded by a polynomial in $n$, $\mathcal{A}$ is called a VC class or a polynomial class of sets. The VC dimension of $\mathcal{A}$, $V_\mathcal{A}$ is the largest positive integer $m$ such that $\Delta_\mathcal{A}(m)=2^m$. If no such $m$ exists then $V_\mathcal{A}= \infty$. The VC density of $\mathcal{A}$, $dV_\mathcal{A}$ is the infimum over the set of positive reals $r$ such that a constant $C$ exists with the property that $\Delta_\mathcal{A}(n) \leq C n^r$ for all $n\in \mathbb{N}$.\\
\end{definition}

The following is one of the cornerstones of the so-called VC theory:
\begin{theorem} \label{vcteo1}(Vapnik-Chervonenkis).
$\mathcal{A}$ is a P-Glivenko-Cantelli class for any measure $P$ if and only if its VC dimension is finite.
\end{theorem}

The if part of proof of \ref{vcteo1} is based on the next results:

\begin{lemma}{(Symmetrization Lemma)}\label{symlem} Let $Z_1, \ldots, Z_n$ and $Z'_1, \ldots, Z'_n$ be i.i.d. random variables taking values in $\mathbb{R}^d$. Denote by $P'_n$ the empirical measure associated to the sample $Z'_1, \ldots, Z'_n$ and by $P_n$ the empirical measure associated to the sample $Z_1, \ldots, Z_n$.
Then for $n\epsilon^2 \geq 2 $ we have
\[
P \left \{ \sup_{A \in \mathcal{A}}|P_n(A)-P(A)| > \epsilon  \right \} \leq 2 P \left \{  \sup_{A \in \mathcal{A}} |P_n(A)-P'_n(A)| > \epsilon/2 \right \}
\]
\begin{proof}
Let $A^* \in \mathcal{A}$ be a set such that $|P_n(A^*)-P(A^*)| > \epsilon$, otherwise let $A^*$ be any fixed set in $\mathcal{A}$. Then
\begin{align*}
P \left \{ \sup_{A \in \mathcal{A}} | P_n(A)-P'_n(A)| > \epsilon/2 \right \} &\geq P \left \{ |P_n(A^*)-P'_n(A^*)|> \epsilon/2 \right \} \\ & \geq P \left \{ | P_n(A^*)-P(A^*)| > \epsilon, |P'_n(A^*) - P(A^*) | < \epsilon/2 \right \} \\ &= P \left \{  1_B P \left \{ | P'_n(A^*)-P(A^*) | < \epsilon/2\left. \right |  Z_1, \ldots, Z_n    \right \}  \right \}
\end{align*}
Where $B$ is the event $ | P_n(A^*)-P(A^*)| > \epsilon $. Bound the conditional probability inside using Chebyshev's inequality
\begin{align*}
P \left \{ |P'_n(A^*)-P(A^*)| < \epsilon/2 \left. \right |  Z_1, \ldots, Z_n \right \} &\geq 1-\dfrac{P(A^*)(1-P(A^*))}{n\epsilon^2/4} \\ & \geq 1- \dfrac{1}{n\epsilon^2} \geq \dfrac{1}{2}
\end{align*}
Whenever $n\epsilon^2 \geq 2$. Then we have
\begin{align*}
P \left \{ \sup_{A \in \mathcal{A}} | P_n(A) - P'_n(A)| > \epsilon/2 \right \} & \geq \dfrac{1}{2} P \left \{ |P_n(A^*)-P(A^*)| > \epsilon \right \} \\ & \geq \dfrac{1}{2} P \left \{ \sup_{A \in  \mathcal{A}} |P_n(A)-P(A)|> \epsilon \right \}
\end{align*}

\end{proof}
\end{lemma}

\begin{lemma}{(Second Symmetrization Lemma)}\label{symlem2} Let $Z_1, \ldots, Z_n$ and $Z'_1, \ldots, Z'_n$ be i.i.d. random variables taking values in $\mathbb{R}^d$. Denote by $P'_n$ the empirical measure associated to the sample $Z'_1, \ldots, Z'_n$. Let $\sigma_1, \ldots, \sigma_n$ be i.i.d. sign variables, independent of the $Z_i$'s and $Z'_i$'s with $P \{ \sigma_i = -1 \} = P \{ \sigma_i=1 \} = 1/2 $, then 
\begin{align*}
P \left \{  \sup_{A \in \mathcal{A}} |P_n(A)-P'_n(A)| > \epsilon/2 \right \} \leq 2 P \left \{ \sup_{A \in \mathcal{A}} \dfrac{1}{n} \left | \sum_{i=1}^n \sigma_i 1_A(Z_i) \right | > \epsilon/4 \right \}
\end{align*}
\begin{proof}

By independence, the next two random variables have the same distribution

\begin{gather*}
\sup_{A\in \mathcal{A}} \left | \sum_{i=1}^n (1_A(Z_i)-1_A(Z'_i)) \right | \\
\sup_{A\in \mathcal{A}} \left | \sum_{i=1}^n \sigma_i (1_A(Z_i)-1_A(Z'_i)) \right |
\end{gather*}
then we have that 
\begin{align*}
P \left \{  \sup_{A \in \mathcal{A}} |P_n(A)-P'_n(A)| > \epsilon/2 \right \} &= P \left \{\sup_{A\in \mathcal{A}} \left | \sum_{i=1}^n (1_A(Z_i)-1_A(Z'_i)) \right | > \epsilon/2 \right \} \\ &= P \left \{\sup_{A\in \mathcal{A}} \left | \sum_{i=1}^n \sigma_i(1_A(Z_i)-1_A(Z'_i)) \right | > \epsilon/2 \right \} \\ & \leq P \left \{\sup_{A\in \mathcal{A}} \left | \sum_{i=1}^n \sigma_i1_A(Z_i) \right | > \epsilon/4 \right \} +  P \left \{\sup_{A\in \mathcal{A}} \left | \sum_{i=1}^n \sigma_i1_A(Z'_i) \right | > \epsilon/4 \right \} \\ &= 2P  \left \{\sup_{A\in \mathcal{A}} \left | \sum_{i=1}^n \sigma_i1_A(Z_i) \right | > \epsilon/4 \right \}
\end{align*}

\end{proof}
\end{lemma}

\begin{lemma}{(Hoeffding's Inequality, for a proof see \cite{pollard84}, Appendix B)}\label{hoefineq} Let $Y_1, \ldots, Y_n$ be independent random variables with zero means and bounded ranges $a_i \leq Y_i \leq b_i$. For each $\eta >0$
$$
P \left \{ | Y_1 + \ldots + Y_n | \geq \eta \right \} \leq 2 \exp \left (  -2\eta^2 / \sum_{i=1}^n (b_i-a_i)^2  \right )
$$ 

\end{lemma}

We have the following bound, which is obtained using simple arguments involving the well known symmetrization lemma and Hoefding's inequality: 
\begin{theorem}{(\cite{devroye96} Theorem 12.5)} \label{edson} 
\[
P \left \{ \sup_{A \in \mathcal{A}} |P_n(A)-P(A)| > \epsilon \right \} \leq 8 \Delta_{\mathcal{A}}(n) \exp(-n\epsilon^2/32)
\]
\begin{proof} By lemmas \ref{symlem} and \ref{symlem2} we have
\begin{align}\label{boundp1}
P \left \{ \sup_{A \in \mathcal{A}} |P_n(A)-P(A)| > \epsilon \right \} \leq 4 P  \left \{\sup_{A\in \mathcal{A}} \left | \sum_{i=1}^n \sigma_i1_A(Z_i) \right | > \epsilon/4 \right \}
\end{align}
To bound the right hand side, condition on the value of $Z_1, \ldots, Z_n$. Upon conditioning, the value of $ \left | \sum_{i=1}^n \sigma_i1_A(Z_i) \right |$ depends on the set $A \cap \{Z_i : i=1, \ldots, n \}$ and there are at most $\Delta_{\mathcal{A}}(n)$ such different sets, by the union bound this implies that
\begin{align*}
P \left \{ \left. \sup_{A \in \mathcal{A}} \dfrac{1}{n} \left | \sum_{i=1}^n \sigma_i 1_A(Z_i) \right | > \epsilon/4 \right | Z_1, \ldots, Z_n \right \} \leq \Delta_{\mathcal{A}}(n) \sup_{A \in \mathcal{A}} P \left \{ \left. \dfrac{1}{n} \left | \sum_{i=1}^n \sigma_i 1_A(Z_i) \right | > \epsilon/4  \right | Z_1, \ldots, Z_n \right \}
\end{align*} 
Given the value of $Z_i$, the variables $\sigma_i 1_A(Z_i)$are independent and take values in $[-1,1]$. Using lemma \ref{hoefineq} we get
\begin{align*}
P \left \{ \left. \dfrac{1}{n} \left | \sum_{i=1}^n \sigma_i 1_A(Z_i) \right | > \epsilon/4  \right | Z_1, \ldots, Z_n \right \} \leq 2 \exp(-n\epsilon^2/32)
\end{align*}
Take expectation and plug into equation (\ref{boundp1}).

\end{proof}
\end{theorem}

\begin{corollary} If $\Delta_{\mathcal{A}}(n)$ is bounded by a polynomial in $n$, then the class $\mathcal{A}$ is $P$-Glivenko-Cantelli for any probability measure $P$.
\begin{proof} Under this condition, the bound of theorem \ref{edson} is summable in $n$. Use Theorem \ref{BC}.

\end{proof}

\end{corollary}

and we have the following bounds for $\Delta_{\mathcal{A}}(n) $ in terms of the VC dimension of the class $\mathcal{A}$: 
\begin{theorem}{(\cite{devroye96} Theorem 13.3)}\label{vcdim} for all $n>2V_\mathcal{A}$
$$
\Delta_{\mathcal{A}}(n) \leq \sum_{i=0}^{V_\mathcal{A}} \binom{n}{i} \leq \left( \dfrac{en}{V_\mathcal{A}} \right )^{V_{\mathcal{A}}}
$$
If $V_\mathcal{A} > 2$ then $\Delta_{\mathcal{A}}(n) \leq n^{V_\mathcal{A}}$. In this case, together with \hyperref[edson]{theorem \ref*{edson}} we have:
\begin{equation}
P \left \{ \sup_{A \in \mathcal{A}} |P_n(A)-P(A)| > \epsilon \right \} \leq 8  n^{V_\mathcal{A}}  \exp(-n\epsilon^2/32)
\end{equation}
\end{theorem}

\begin{remark} The authors of \cite{adolfo2014} point out that the actual quantity involved in the proof of such theorem is the VC density of the class, hence, a class is P-Glivenko-Cantelli for any probability measure P if and only if its VC density is finite. Also, bounds in terms of the $VC$ dimension of a class, can be replaced by analogous bounds in terms of the VC density. For example, theorem \ref{edson} becomes the following:

$$
P \left \{ \sup_{A \in \mathcal{A}} |P_n(A)-P(A)| > \epsilon \right \} \leq 8 C n^{dV_\mathcal{A}}\exp(-n\epsilon^2/32)
$$

where the constant $C$ depends on the class $\mathcal{A}$. 
\end{remark}
\end{subsection}

\begin{subsection}{Empirical Risk Minimization}\label{Empirical}
Let $X$ be a random variable taking values on a set $S \subset \mathbb{R}^d$ and $Y$ a random variable taking values on the set $\{0,1\}$, both defined on the same probability space. Let $\mathcal{G}$ be a class of functions with domain $S$ taking values in $\{0,1\}$ and let $P$ be the joint distribution of $(X,Y)$.

\begin{definition}\label{risks}
Define the risk of a function $g\in \mathcal{G}$ as 
\begin{align*}
L(g)=P\left (|g(X)-Y| \right)
\end{align*}
Given an i.i.d. random sample $(X_i,Y_i), \ldots, (X_n,Y_n)$ from the distribution of $(X,Y)$, define the empirical risk of a function $g \in \G$:
\begin{align*}
\widehat{L}_n(g)=\dfrac{1}{n} \sum_{i=1}^n |g(X_i)-Y_i|
\end{align*}
\end{definition}

For a class $\G$ we would like to bound the quantity $P \left \{ \sup_{g \in \mathcal{G}} |\widehat{L}_n(g)-L(g)| > \epsilon \right \}$. This measures how far is the empirical risk of a function $g$, directly computable from a random sample $(X_i,Y_i)$, from the actual risk, which is not directly computable in most practical situations.

\begin{definition}\label{class} Let $\G$ be a class of functions of the form $g:S \rightarrow \{0,1\}$. Define $\tilde{\mathcal{A}}$ as the collection of sets
\begin{align*}
\{ \{ x : g(x)=1\} \times \{0\}\} \bigcup \{ \{ x : g(x)=0\} \times \{1\}\}, \, \, \, g \in  \G
\end{align*} 
\end{definition}

\begin{lemma}\label{relation}
\begin{align*}
\sup_{g \in \G} | \widehat{L}_n(g)-L(g)| = \sup_{A \in \tilde{A}} |P_n(A)-P(A)|
\end{align*}
\begin{proof}

\begin{gather*}
L(g)=P\left ( |g(X)-Y| \right)=P(g(X) \neq Y )= P\left(\{ x : g(x) =1\} \times \{0\} \bigcup \{ x : g(x) =0\} \times \{1\} \right ) \\
\widehat{L}_n(g)= \dfrac{1}{n} \sum_{i=1}^n  |g(X_i)-Y_i| =  \dfrac{1}{n} \sum_{i=1}^n 1_{\{ g(X_i) \neq Y_i\}} = P_n \left(\{ x : g(x) =1\} \times \{0\} \bigcup \{ x : g(x) =0\} \times \{1\} \right )
\end{gather*}
\end{proof}
\end{lemma}

\begin{theorem}{(\cite{devroye96} Theorem 12.6)}\label{riskbound} For a  class of functions $\mathcal{G}$ with domain $S$ taking values in $\{0,1\}$, and $\mathcal{A}$ defined as before
$$
P \left \{ \sup_{g \in \mathcal{G}} |\widehat{L}_n(g)-L(g)| > \epsilon \right \} \leq 8 \Delta_{\tilde{\mathcal{A}}} (n) \exp(-n\epsilon^2/32)
$$

In the case where $2<V_{\tilde{\mathcal{A}}}<\infty$ and $n> 2V_{\tilde{\mathcal{A}}}$, by \hyperref[vcdim]{theorem \ref*{vcdim}}

$$
P \left \{ \sup_{g \in \mathcal{G}} |\widehat{L}_n(g)-L(g)| > \epsilon \right \} \leq 8 n^{V_{\tilde{\mathcal{A}}}} \exp(-n\epsilon^2/32)
$$
\begin{proof}
Direct consequence of lemma \ref{relation} and theorems \ref{edson} and \ref{vcdim}.

\end{proof}

\end{theorem}

\begin{lemma}{(For a proof see \cite{devroye96}, Theorem 13.1.)} \label{simp}
Let $ \mathcal{A}$ be the class of sets of the form $\{x : g(x)=1 \}$ with $g  \in \mathcal{G}$ then for every $n$, $\Delta_{\tilde{\mathcal{A}}}(n) = \Delta_{\mathcal{A}}(n)$ and hence $V_{\tilde{\mathcal{A}}}=V_{\mathcal{A}}$.

\end{lemma}

\begin{definition} $V_{\G}= V_{\tilde{\mathcal{A}}}= V_{\mathcal{A}}$.

\end{definition}

In section \ref{Means} we will further introduce analogous bound in the context of Regression Estimation where the variable $Y$ takes values in the interval $[0,1]$.
\end{subsection}

\begin{subsection}{Vapnik's $\eta$ trick}\label{Trick}
\begin{definition} Let $\hat{g}_n$ be a classifier in $\mathcal{G}$ minimizing the empirical risk over the class.
\begin{gather*}
\hat{g}_n = \argmin_{\mathcal{G}} \widehat{L}_n(g)
\end{gather*}
\end{definition}

We will obtain a bound for the real risk of this classifier in term of its empirical risk and the VC dimension of the class. 

\begin{theorem} with probability at least $1-\eta$ simultaneously for all $g \in \mathcal{G}$
\[
L(g) \leq \widehat{L}_n(g)+ \sqrt{\dfrac{32(V_\mathcal{G}\log(n)-\log(\eta/8))}{n}}
\]
in particular for $\hat{g}_n$
\begin{equation}\label{stbound}
L(\hat{g}_n)\leq \widehat{L}_n(\hat{g}_n)+ \sqrt{\dfrac{32(V_\mathcal{G}\log(n)-\log(\eta/8))}{n}}
\end{equation}
\begin{proof}
Let $\eta=8 n^{V_\mathcal{G}} \exp(-n\epsilon^2/32)$. Solving for $\epsilon$ we get
\[
\epsilon= \sqrt{\dfrac{32(V_\mathcal{G}\log(n)-\log(\eta/8))}{n}}
\]
Then by \hyperref[riskbound]{theorem \ref*{riskbound}} we have that with probability at least $1-\eta$ and any $g \in \G$
\[
 L(g)-\widehat{L}_n(g)\leq |L(g)-\widehat{L}_n(g)| \leq \sup_{g \in \mathcal{G}} |L(g)-\widehat{L}_n(g)| \leq \sqrt{\dfrac{32(V_\mathcal{G}\log(n)-\log(\eta/8))}{n}}
\]

\end{proof}

\end{theorem}

We also have the following 
\begin{lemma} \label{Lbound} Let $\mathcal{C}$ be a class of $\{0,1\} $ mappings, $(X_1,Y_1), \ldots, (X_n,Y_n)$ a random sample from the distribution of $(X,Y)$ and $\phi^*_n$ a classifier minimizing the empirical risk over the class, then
\[
L(\phi^*_n) - \inf_{\phi \in \mathcal{C}} L(\phi) \leq 2 \sup_{\phi \in \mathcal{C}}| \widehat{ L}_n(\phi)-L(\phi)|
\]
\begin{proof}
\begin{align*}
L(\phi^*_n)-\inf_{\phi \in \mathcal{C}} L(\phi) &= L(\phi^*_n)-\widehat{L}_n(\phi^*_n)+\widehat{L}_n (\phi^*_n)-\inf_{\phi \in \mathcal{C}} L(\phi) \\
 & \leq L(\phi^*_n) - \widehat{L}_n(\phi^*_n) + \sup_{\phi \in \mathcal{C}} | \widehat{L}_n(\phi) - L(\phi)| \\
 & \leq 2 \sup_{\phi \in \mathcal{C}}|\widehat{L}_n(\phi) - L (\phi)|
\end{align*}
\end{proof}
\end{lemma}

\begin{theorem}{(\cite{devroye96} Theorem 12.6)}\label{riskbound2}
\begin{gather*}
P \left \{ L(\hat{g}_n)-\inf_{g\in \mathcal{G}} L(g) > \epsilon \right \} \leq 8 n^{V_\mathcal{G}} \exp(-n\epsilon^2/128)
\end{gather*}
and with probability at least $1-\eta$

\begin{equation}
L(\hat{g}_n) \leq \inf_{g\in \mathcal{G}} L(g) +  \sqrt{\dfrac{128(V_\mathcal{G}\log(n)-\log(\eta/8))}{n}}
\end{equation}
\end{theorem}

In sections \hyperref[VC_Subgraph]{\ref*{VC_Subgraph}} and \hyperref[Parametric]{\ref*{Parametric}} we will generalize this method to obtain such type of bounds in the context of Regression Estimation.

\end{subsection}
\begin{subsection}{Structural Risk Minimization}\label{Structural}
Let $D_n=\{(X_1,Y_1), \ldots, (X_n,Y_n)\}$ be a random sample of size $n$ from the distribution of a random pair $(X,Y)$ taking values in $S \times \{0,1\}$. A map $D_n \mapsto g_n$ is called a \emph{classification rule} and $L(g_n)$ is denoted as $L_n$.
\begin{definition}{(Bayes Risk)}
For a given distribution of $(X,Y)$ define the Bayes regression function as a function $g^*: S \rightarrow \{0,1\}$ such that $L(g^*) \leq L(g)$ for any measurable function $g$. Define $L^*=L(g^*)$ as the Bayes Risk.
\end{definition}

In practice, as the distribution of $(X,Y) $ is unknown, so is the Bayes classifier $g^*$, which means we can only aim to find a function $g_n$ based on the sample that is close to $g^*$ in some way. This motivates the following
\begin{definition}\label{consistency} A regression rule is consistent for a certain distribution of $(X,Y)$ if 
$$
E(L_n) \rightarrow L^*,  \, \, \, \, n \rightarrow \infty
$$
It is strongly consistent if
$$
\lim_{n\to \infty} L_n = L^* \, \, \, \, a.s.
$$
and it is called strongly universally consistent, if it is strongly consistent for any distribution of $(X,Y)$.
\end{definition}

Let $\mathcal{G}_1, \mathcal{G}_2, \ldots$ be a sequence of classes of classifiers. Let $\hat{g}_{n,j}$ be a classifier in $\mathcal{G}_j$ minimizing the empirical error $\widehat{L}_n$. Let $g^*_n$ be the classifier among $\{ \hat{g}_{n,j} : j =1,2, \ldots \}$ minimizing the term
$$
\widetilde{L}_{n,j}= \widehat{L}_n(\hat{g}_{n,j})+ \sqrt{\dfrac{32V_{\mathcal{G}_j}\log(en)}{n}}
$$
This is called the classifier based on \emph{Structural Risk Minimization}.

\begin{theorem}\label{strong consistency}\hyperref[devroye96]{(\cite{devroye96}, Theorem 18.2)} Suppose that for any distribution of $(X,Y)$
$$
\lim_{j \to \infty} \inf_{g \in \mathcal{G}_j} L(g) = L^*
$$
Suppose also that the VC dimensions of $\mathcal{G}_1, \mathcal{G}_2, \ldots$ are finite and satisfy
$$
\Delta = \sum_{j=1}^\infty e^{-V_{\mathcal{G}_j}} < \infty
$$
Then the classification rule $g^*_n$ based on structural risk minimization, is strongly universally consistent.
\end{theorem}
In Chapter \ref{Strong} we will state and prove a more general result based on the previous theorem. The idea behind this approach is to minimize the sum of the empirical error on the sample and a complexity term $r(j,n)$ depending on the size of the sample and the index $j$ of the class. With increasing $j$, each class $\mathcal{G}_j$ is a more complex class with larger VC dimension, leading to a potential decrease in empirical error over the sample. However choosing a class with complexity too large, may lead to a higher risk, this is known as the problem of \emph{overfitting}. The theorem shows that the additional complexity term prevents this from happening. This comes from equation \hyperref[stbound]{(\ref*{stbound})} where we got that the empirical error plus an $O(\sqrt{V_\mathcal{G}\log(n)/n})$ term bounds the real risk. Adding a term of the same magnitude prevents from underestimating such risk.\\

 In section \hyperref[Strong]{\ref*{Strong}} we will prove a similar result in the context of Regression Estimation.
 \end{subsection}
\end{section}

\begin{section}{Regression Estimation}\label{Regression}
\vspace{5 mm}

Just as P-Glivenko-Cantelli classes of sets are the core of classifier selection, P-Glivenko-Cantelli classes of functions are the basis of regression function estimation. However a characterization of such classes is far more difficult than in the case of classes of sets. We will state some sufficient conditions based on concepts from metric spaces. All the results and definitions also work in the case where we have a pseudometric space. We might refer to a pseudometric space as a metric space if no confusion arises.\\

\begin{subsection}{Metric Space Theory}\label{Metric_Space}
\begin{definition}{(Covering Number)}\label{coveringnumber}  Let $(M,d)$ be a metric space, $\epsilon > 0$. A set $B\subset M$ is an $\epsilon$-cover for $M$ if, for every $x \in M$, there exists some $b \in B$ such that $d(x,b) < \epsilon$. The $\epsilon$-covering number of M, $N(\epsilon,M,d)$ is the minimal cardinality of an $\epsilon$-cover for $M$. It is defined to be $\infty$ if no such finite cover exists. Let $P$ be a probability measure on a set $S$ and $\mathcal{F}$ be a set of $P$-measurable functions with domain $S$. Define $N_p(\epsilon, \mathcal{F},P)$ to be the $\epsilon$-covering number of $\mathcal{F}$ with respect to the $L^p(P)$ metric and $\log N_p(\epsilon, \mathcal{F}, P)$ to be the \emph{$\epsilon$-metric entropy} of the class.
\end{definition}

\begin{remark}\label{remcov} The $\epsilon$-covering number of a subset $A\subset M$ of a metric space is defined as the minimum cardinality of a set $\{x_i \}$ of elements of $M$ such that the balls with center $x_i$ and radius $\epsilon$ cover $A$. Denote such covering number as $N(\epsilon, A\subset M, d)$. Regarding the $\epsilon$ covering numbers of $A$ as a metric space on its own we get the following inequality:
$$
N(\epsilon, A, d) \leq N(\epsilon/2, A\subset M, d)
$$ 
\end{remark}

\begin{definition}Let Let $(M_1,d_1)$ and $(M_2,d_2)$ be two metric spaces and $f:M_2 \rightarrow M_1$ a function such that there exist constants $\alpha > 0$, $K > 0 $ such that 
$$
d_1(f(x), f(y)) \leq K d_2(x,y)^\alpha
$$
$f$ is called H\"older continuos of order $\alpha$ with constant $K$. If $\alpha=1$ the function is called Lipschitz continuos.
\end{definition}

\begin{lemma} \label{metrics1}Let $(M_1,d_1)$ and $(M_2,d_2)$ be two bounded metric spaces and a surjective H\"older continuous function of order $\alpha$ with constant $K$ $f:M_2 \rightarrow M_1$, then
$$
N(K\epsilon^\alpha, M_1, d_1) \leq N(\epsilon, M_2, d_2)
$$
\begin{proof}
let $\{x_1, \ldots x_m\} \subset M_2$ be a set of minimum cardinality such that for any $x \in M_2$ there exists $x_j$ such that $d_2(x,x_j) < \epsilon$, then $d_1(f(x),f(x_j)) \leq K d_2(x,x_j)^\alpha < K\epsilon^\alpha$. This implies that $N(K\epsilon^\alpha, d_1, M_1) \leq |f(x_1), \ldots, f(x_m)| \leq m =N(\epsilon, d_2, M_2)$.
\end{proof}
\end{lemma}

\begin{definition} A function $G$ is called an envelope of a class of functions $\mathcal{G}$ if $|g(x)| \leq G(x)$, $\forall f \in \mathcal{G}$.
\end{definition}
For ease of computations we will assume the class $\mathcal{G}$ comprises of positive functions bounded by 1, so that the constant function $1$ is an envelope for the class. Results easily generalize for classes of functions with constant envelope. 

\begin{lemma}\label{squared1} Let $\mathcal{G}$ be a class of positive functions with domain $S$ and constant envelope $1$. Then for any probability measure $P$ on $S$
$$
N_1(2\epsilon, \mathcal{G}^2, P) \leq N_2(\epsilon, \mathcal{G}, P)
$$
where $\mathcal{G}^2=\{g^2: g \in \mathcal{G}\}$ is the class of squared functions.
\begin{proof}
\[
P\left (|g_1^2-g_2^2| \right )=P\left( |(g_1-g_2)(g_1+g_2)| \right ) \leq \sqrt{P(|g_1-g_2|)^2}\sqrt{P(|g_1+g_2|)^2} \leq 2 \sqrt{P(|g_1-g_2|)^2}
\]

by virtue of the Cauchy-Schwarz inequality. Use lemma \hyperref[metrics1]{(\ref*{metrics1})} for the function $g \mapsto g^2$ which is Lipschitz continuous with constant 2.

\end{proof}
\end{lemma}
This means that the $N_1$ entropy of the square of a class of functions with constant envelope can be bounded by the $N_2$ entropy of the class.
\begin{lemma}
Let $\mathcal{G}$ be a class of functions with domain $S$ and $h$ a fixed function with the same domain but possibly not in $\mathcal{G}$. Let  $\mathcal{G}+h=\{g+h: g \in \mathcal{G} \}$. Then for any probability measure P, if $h \in  L_p(P)$
$$
N_p(\epsilon, \mathcal{G}, P)=N_p(\epsilon, \mathcal{G}+h,P)
$$
\begin{proof} Suppose that $\{g_1, \ldots, g_m \}$ is minimal $\epsilon$-cover for $\G$ with metric $L_p(P)$. Let $g+h$ be any function in $\G+h$. $\| g-g_i \|_p < \epsilon$ for some $g_i$. Then $\|g+h-(g_i+h)\|_p = \|g-g_i\|_p < \epsilon$. This shows that the set $\{g_1+h, \ldots, g_m+h\}$ is an $\epsilon$-cover for the class $\G+h$ with metric $L_p(P)$ so that $N_p(\epsilon, \mathcal{G}+h,P) \leq N_p(\epsilon, \mathcal{G}, P)$. A symmetric argument yields the result.

\end{proof}
\end{lemma}

\begin{lemma}\label{preservemetric} Let $(M_1,d_1)$ and $(M_2, d_2)$ be two (pseudo)metric spaces, and let $f:M_1\rightarrow M_2$ be a bijection. Suppose that $f$ is Bilipschitz with constant $K \geq 1$, that is
$$
\dfrac{1}{K}d_1(x,y) \leq d_2(f(x),f(y))\leq  K d_1 (x,y)
$$
for all $x,y \in M_1$. Then 
$$
N(K\epsilon, M_2, d_2) \leq N(\epsilon, M_1, d_1) \leq N(\epsilon/K, M_2, d_2)
$$
In particular if $K=1$,  $N(\epsilon, M_1, d_1)=N(\epsilon, M_2, d_2)$.
\begin{proof}
We have that $d_2(f(x),f(y)) \leq Kd_1(x,y)$, as $f$ is surjective, we have by lemma \ref{metrics1} that $N(K\epsilon, M_2,d_2) \leq N(\epsilon, M_1, d_1)$. On the other hand, $d_1(x,y) \leq K d_2(f(x),f(y))$. let $a=f(x)$ and $b=f(y)$ then $d_1(f^{-1}(a), f^{-1}(b)) \leq K d_2(a.b)$. As $f^{-1}$ is surjective, we have by lemma \ref{metrics1} that $N(K\epsilon, M_1, d_1) \leq N(\epsilon, M_2, d_2)$ so that $N(\epsilon, M_1, d_1) \leq N(\epsilon/K, M_2, d_2)$.

\end{proof}
\end{lemma}

\begin{definition} Let $(M,d)$ be a metric space and $x\in M$.
\[
B(\delta, x, d)= \{ y \in M: d(y,x) < \delta \}
\]

\end{definition}

\begin{lemma}\label{metric1}
Let $(M_1,d_1)$ and $(M_2, d_2)$ be two metric spaces. $f:M_1\rightarrow M_2$ be a surjective map and $x_0 \in M_1$ a fixed element. Suppose that for any $x \in M_1$ 
$$
d_1^2(x, x_0) \leq  d_2(f(x),f(x_0))
$$
Then $B(\delta, f(x_0), d_2) \subset f(B(\sqrt{\delta}, x_0, d_1))$.
\begin{proof}
Let $f(x) \in M_2$ such that $d_2(f(x),f(x_0))< \delta$,  then $d_1^2(x,x_0) \leq d_2(f(x),f(x_0)) < \delta $. This implies that $d_1(x,x_0) < \sqrt{\delta}$. So $x \in B(\sqrt{\delta}, x_0, d_1)$.

\end{proof}

\end{lemma}

\end{subsection}

\begin{subsection}{Uniform Convergence of Means to Expectations}\label{Means}
Focusing on the case of a class of positive functions bounded by 1, we have the following theorem, which is a direct adaptation of \hyperref[pollard84]{\cite{pollard84}} Chapter 2, Theorem 24.

\begin{theorem} \label{bound} Let $\mathcal{F}$ be a permissible class of positive functions bounded by 1.
\[
P \{||P_n-P|| > \epsilon \} \leq 8 P \left ( N_1 (\epsilon/8, \mathcal{F}, P_n) \right ) \exp(-n\epsilon^2/128)
\]
where we denote
 \[
 ||P_n-P||= \sup_{f \in \mathcal{F}} |P_n(f)-P(f)|
 \]
 A sufficient condition for a law of large numbers to hold over $\mathcal{F}$ is that $\log N_1 (\epsilon, \mathcal{F}, P_n)/n \rightarrow 0 $ in probability as $n\rightarrow \infty$, for any $\epsilon>0$ .

\end{theorem}

For the proof of this theorem we need the following lemmas which are a direct generalization of lemmas \ref{symlem} and \ref{symlem2}.

\begin{lemma}{(Symmetrization lemma)}\label{symlemf}  Let $Z_1, \ldots, Z_n$ and $Z'_1, \ldots, Z'_n$ be i.i.d. random variables taking values in $\mathbb{R}^d$. Denote by $P'_n$ the empirical integral associated to the sample $Z'_1, \ldots, Z'_n$ and $P_n$ the empirical integral associated to the sample $Z_1, \ldots, Z_n$.
Then for $n\epsilon^2 \geq 2 $ we have
\[
P \left \{ \sup_{f \in \mathcal{F}}|P_n(f)-P(f)| > \epsilon  \right \} \leq 2 P \left \{  \sup_{f \in \mathcal{F}} |P_n(f)-P'_n(f)| > \epsilon/2 \right \}
\]
\begin{proof} The proof of lemma \ref{symlem} carries verbatim replacing $A^*$ by a function $f^*$ such that $|P_n(f^*)-P(f^*)|>\epsilon$, replacing $A \in \mathcal{A}$ by $f \in \mathcal{F}$ and noting that 
\[ 
\var(P'_n(f^*)) =  \dfrac{1}{n} P(f^*)(1-P(f^*))
\]
Wich has maximum value $\frac{1}{4n}$ as $0 \leq P(f^*) \leq 1$.

\end{proof}

\end{lemma}

\begin{lemma}{(Second symmetrization lemma)} \label{symlem2f}
Let $Z_1, \ldots, Z_n$ and $Z'_1, \ldots, Z'_n$ be i.i.d. random variables taking values in $\mathbb{R}^d$. Denote by $P'_n$ the empirical integral associated to the sample $Z'_1, \ldots, Z'_n$ and by $P_n$ the empirical integral associated to the sample $Z_1, \ldots, Z_n$. Let $\sigma_1, \ldots, \sigma_n$ be i.i.d. sign variables, independent of the $Z_i$'s and $Z'_i$'s with $P \{ \sigma_i = -1 \} = P \{ \sigma_i=1 \} = 1/2 $, then 
\begin{align*}
P \left \{  \sup_{f \in \mathcal{F}} |P_n(f)-P'_n(f)| > \epsilon/2 \right \} \leq 2 P \left \{ \sup_{f \in \mathcal{F}} \dfrac{1}{n} \left | \sum_{i=1}^n \sigma_i f(Z_i) \right | > \epsilon/4 \right \}
\end{align*}
\begin{proof} The proof of lemma \ref{symlem2} carries verbatim replacing $A \in \mathcal{A}$ by $f \in \mathcal{F}$ and $1_A$ by $f$.

\end{proof}

\end{lemma}

\textbf{Proof of theorem \ref{bound}:}
\begin{proof}
The two lemmas \ref{symlemf} and \ref{symlem2f} yield, for $n\epsilon^2 \geq 2$
\begin{equation}\label{prebound}
P \{||P_n-P|| > \epsilon \} \leq 4 P \left \{ \sup_{f \in \mathcal{F}} \dfrac{1}{n} \left | \sum_{i=1}^n \sigma_i f(Z_i) \right | > \epsilon/4 \right \}
\end{equation}
Now given $Z_1, \ldots, Z_n$ let $\{ f_1, \ldots, f_m \}$, where $m=N_1(\epsilon/8, \mathcal{F}, P_n)$, such that for any $f \in \mathcal{F}$, $P_n(|f-f_j|) < \epsilon/8$ for some $j$. Write $f^*$ for such $f_j$. Write $P^\circ_n(f)=  \frac{1}{n} \sum_{i=1}^n \sigma_i f(Z_i) $. Note that
\[
|P^\circ_n(f)| =   \left| \dfrac{1}{n} \sum_{i=1}^n \sigma_i f(Z_i) \right|  \leq \dfrac{1}{n} \sum_{i=1}^n |f(Z_i)| = P_n (|f|)
\]
Using this and working out the right hand side of equation (\ref{prebound}) we have
\begin{align*}
P \left \{  \left. \sup_{f \in \mathcal{F}} | P^\circ_n(f) | > \epsilon/4 \right | Z_1, \ldots, Z_n \right \} &\leq P \left \{ \left. \sup_{f \in \mathcal{F}} |P^\circ_n(f^*+f-f^*)| > \epsilon/4 \right | Z_1, \ldots, Z_n \right \} \\ & \leq P \left \{ \left. \sup_{f \in \mathcal{F}} |P^\circ_n(f^*)|+P_n|f-f^*| > \epsilon/4 \right | Z_1, \ldots, Z_n \right \} \\ & \leq P \left \{ \left. \sup_{f \in \mathcal{F}} |P^\circ_n(f^*)| > \epsilon/8 \right | Z_1, \ldots, Z_n \right \} \\ & \leq P \left \{ \left. \max_{j} |P^\circ_n(f_j)| > \epsilon/8 \right | Z_1, \ldots, Z_n \right \}  \\ & \leq N_1(\epsilon/8, \mathcal{F}, P_n) \max_{j} P \left \{ \left. |P^\circ_n(f_j)| > \epsilon/8 \right | Z_1, \ldots, Z_n \right \} 
\end{align*}
Finally, using again Hoeffding's inequality we have that
\begin{align*}
P \left \{ \left. |P^\circ_n f_j | > \epsilon/8 \right | Z_1, \ldots, Z_n \right \} &= P \left \{ \left. \left | \sum_{i=1}^n \sigma_i f_j(Z_i) \right| > n\epsilon/8 \right | Z_1, \ldots, Z_n \right \} \\ & \leq 2 \exp \left (  -2 (n\epsilon/8)^2 / \sum_{i=1}^n (2f_j(Z_i) )^2 \right ) \\ &\leq 2 \exp \left ( -n \epsilon^2 / 128\right)
\end{align*}
Take expectation to get rid of the conditional and plug into equation (\ref{prebound}).
\end{proof}

\begin{definition}
In the case of regression estimation, given a random sample $(X_1,Y_1),\ldots, (X_n,Y_n)$ where $Y$ is a $[0,1]$ valued random variable, $X$ takes values on a set $S$, and $\mathcal{G}$ is a class of positive functions bounded by 1 with domain $S$, we will define our risk as the expected squared error
\[
L(g)= P\left((g(X)-Y)^2\right )
\]
and the empirical risk will be
\[
\widehat{L}_n(g)=\dfrac{1}{n} \sum_{i=1}^n (g(X_i)-Y_i)^2 = P_n((g(X)-Y)^2)
\]
The mean square error over the sample. 
\end{definition}

\begin{definition}\label{loss class} Let $X$ be a random variable taking values in a set $S\subset \mathbb{R}^d$ and $Y$ be a random variable taking values in $[0,1]$. For a class of functions $\G$ with domain $S$ and range $[0,1]$ define, for $g \in  \G$, its loss function $ l_g:S \times  [0,1] \rightarrow \mathbb{R}$, $(x,y) \mapsto \left(g(x)-y\right )^2$, and the loss class
\[
\mathcal{L}_{\G}= \{ l_g : g \in \G \}
\]

\end{definition}

\begin{lemma}\label{env1}
Let $(X,Y)$ be a $S\times[0,1]$ valued random pair. $\mathcal{G}$ a class of positive functions with domain $S$ and bounded by $1$. Then the class of functions $\mathcal{L}_\mathcal{G}$ with domain $S\times [0,1]$ is a class of positive functions bounded by $1$.
\begin{proof} As $0\leq g(x) \leq 1$ and $0\leq y \leq 1$, $0 \leq |g(x)-y| \leq 1 $ and $0 \leq (g(x)-y)^2 \leq 1$, for all $x \in S$ and $g \in \mathcal{G}$.

\end{proof}
\end{lemma}

\begin{theorem} \label{importantbound}
\begin{equation}
P \left \{ \sup_{g \in \mathcal{G}} |L(g)-\widehat{L}_n(g)| > \epsilon \right \} \leq 8 P N_1(\epsilon/8, \mathcal{L}_\mathcal{G}, P_n) \exp(-n\epsilon^2/128)
\end{equation}
\begin{proof}
Note that 
\begin{gather*}
L(g)= P(l_g)\\
\widehat{L}_n(g)= P_n(l_g)\\
\sup_{g\in \G} |L(g)-\widehat{L}_n(g)| = \sup_{l_g \in \mathcal{L}_{\G}} |P(l_g)-P_n(l_g)|
\end{gather*}

\end{proof}
Using theorem \hyperref[bound]{(\ref*{bound})} we get the result.
\end{theorem}

\begin{definition} For a class $\G$ of positive functions bounded by $1$ with domain $S$, a random pair $(X,Y)$ taking values in $S \times [0,1]$ and an i.i.d. random sample $(X_1,Y_1), \ldots, (X_n,Y_n)$ from the distribution of $(X,Y)$ let $\hat{g}_n$ be a function in $\G$ such that
\[
\widehat{L}_n(g) = \min_{g \in \G} \widehat{L}_n(g)
\]
We assume such function exists. Call this function the least squares estimator.
\end{definition}

We want to obtain a confidence interval for $|L(g)-\widehat{L}_n(g)|$, the deviation of the empirical risk from the real risk, and a confidence interval for $L(\hat{g}_n)- L(g_0)$, the deviation of the real risk of the least squares estimator to that of the best function in the class.

\begin{lemma} \label{Lbound} Let $\mathcal{G}$ be a class of functions taking values in $[0,1]$.
$$
L(\hat{g}_n) - \inf_{g \in \mathcal{G}} L(g) \leq 2 \sup_{g \in \mathcal{G}}| \widehat{ L}_n(g)-L(g)|
$$
\begin{proof}
The same as the proof of lemma \ref{Lbound}.
\end{proof}

\end{lemma}

This lemma shows that bounds for $ \sup_{g \in \mathcal{G}}|\widehat{L}_n(g) - L (g)| $ yield bounds for the efficiency of the least squares estimator. The following corollary makes this precise:
\begin{corollary}
$$
P\left \{ L(\hat{g}_n)-\inf_{g \in \G} L(g)   > \epsilon \right \} \leq P \left \{ \sup_{g \in \mathcal{G}} |L(g)-\widehat{L}_n(g)| > \epsilon/2 \right \} \leq 8 P N_1(\epsilon/16, \mathcal{L}_\mathcal{G}, P_n) \exp(-n\epsilon^2/512)
$$

\end{corollary}

\end{subsection}

\begin{subsection}{VC Subgraph Classes}\label{VC_Subgraph}
There is a large class of function classes for which the $L_r(P)$ $\epsilon$-metric entropy is bounded by a polynomial in $\epsilon^{-1}$ that depends only on the class of functions and $r$, for any probability measure. We refer to the well known VC subgraph classes. We show that if $\mathcal{G}$ is a VC subgraph class, then the class of squared errors $\mathcal{L}_{\G}=\{(g(x)-y)^2: g \in \mathcal{G} \}$ is also a VC subgraph class and hence the $L_1(P)$ $\epsilon$-covering numbers are bounded by a polynomial in $\epsilon^{-1}$.
\begin{definition} Let $\mathcal{F}$ be a class of real valued functions. The set
$$
\subg(f) = \{ (x,t) : t < f(x) \}
$$
is called the subgraph of $f$ and $\subg(\mathcal{F})= \{ \subg{(f)} : f \in \mathcal{F} \}$ is called the subgraph class of $\mathcal{F}$.
\end{definition}

\begin{definition} Call a class of functions $\mathcal{F}$ a VC subgraph class if $\subg(\mathcal{F})$ is a VC class. The VC subgraph dimension of this class $V_\mathcal{F}$ is defined as the VC dimension of its class of subgraphs.
\end{definition}

The following lemma is useful for computing VC subgraph dimensions of classes of functions from the VC subgraph dimension of other classes.
\begin{lemma}\label{subgbounds} (\cite{vaart96}, 2.6.18) Let $\mathcal{F}$ and $\mathcal{G}$ be two VC subgraph classes with domain $S$. $y: S\rightarrow \mathbb{R}$, $\phi: \mathbb{R} \rightarrow \mathbb{R}$ be fixed functions. Denote $\wedge$ and $\vee$ as minimum and maximum, respectively.Then the following hold
\begin{enumerate}
\item $ \mathcal{F} \wedge \mathcal{G} = \{ f \wedge g : f \in \mathcal{F}, g\in \mathcal{G} \}$ is VC subgraph with $V_{\mathcal{F} \wedge \mathcal{G}} \leq V_\mathcal{F}+V_\mathcal{G} -1$.
\item $ \mathcal{F} \vee \mathcal{G} = \{ f \vee g : f \in \mathcal{F}, g\in \mathcal{G} \}$ is VC subgraph with $V_{\mathcal{F} \vee \mathcal{G} } \leq V_\mathcal{F}+V_\mathcal{G} -1$.
\item $-\mathcal{F}$ is VC subgraph with $V_{-\mathcal{F}}=V_\mathcal{F}$.
\item $\mathcal{F}+y=\{f+y: f \in \mathcal{F} \} $ is VC subgraph with $V_{\mathcal{F}+y}=V_\mathcal{F}$.
\item $ \phi \circ \mathcal{F} $ is VC subgraph with $V_{\phi \circ \mathcal{F}} \leq V_\mathcal{F}$ for monotone $\phi$.
\end{enumerate}
\end{lemma}

\begin{lemma}\label{squared} let $\mathcal{G}$ be a VC subgraph class of functions. Then the class of square errors $\mathcal{L}_{\G}$ is a $VC$ subgraph class with $V_{\mathcal{L}_{\G}} \leq 2 V_{\mathcal{G}}-1 $.
\begin{proof} By lemma \ref{subgbounds} (4.) if $\mathcal{G}$ is VC subgraph then $\mathcal{G}-y$ is VC subgraph with the same index, so it suffices to show that $\mathcal{G}^2$ is VC subgraph whenever $\mathcal{G}$ is.\\

Let $\{ (x_1,t_1), \ldots, (x_n, t_n) \} $ be a maximal set of points shattered by $\subg(\mathcal{G}^2)$. That is $n=V_{\mathcal{G}^2}$. The class of squares is a class of positive functions and hence $t_i > 0$ for all $i$, otherwise $(x_i,t_i)$ would lie on the subgraph of all the functions in the class. now if $t_i < g^2(x_i)$ then $\sqrt{t_i} < g(x_i)$ or $\sqrt{t_i} < -g(x_i) $, so $\sqrt{t_i} < g \vee -g (x_i)$ and the set $\{(x_1, \sqrt{t_1}), \ldots, (x_n, \sqrt{t_n})\}$ is shattered by the subgraph of $\{g \vee -g : g \in  \mathcal{G} \} \subset \{ g \vee h: g \in  \mathcal{G}, h \in -\mathcal{G} \} $ then by monotonicity of the subgraph dimension and lemma \ref{subgbounds} (2.) the result follows.

\end{proof}
\end{lemma}

\begin{lemma}\label{abs} Let $\mathcal{G}$ be a VC subgraph class, then $|\mathcal{G}|$ is a $VC$ subgraph class with $V_{|\mathcal{G}|}  \leq 2 V_\mathcal{G}-1 $.
\begin{proof} In lemma \ref{squared} we proved that $V_{\mathcal{G}^2} \leq  2 V_\mathcal{G}-1$ and by monotonicity of the square root and lemma \ref{subgbounds} (5.) we get that
$$
V_{|\mathcal{G}|}= V_{\sqrt{\mathcal{G}^2}} \leq V_{\mathcal{G}^2} \leq  2 V_\mathcal{G}-1
$$
\end{proof}
\end{lemma}

We have the following way to bound the metric entropy depending on the VC subgraph dimension of the class.
\begin{theorem}\label{vcsubgraph}{\hyperref[vaart96](\cite{vaart96}, Theorem 2.6.7)} For a VC class of functions $\mathcal{F}$ with constant envelope 1, for any probability measure P
$$
N_r(\epsilon, \mathcal{F}, P) \leq K (V_\mathcal{F}+1)(16e)^{V_\mathcal{F}+1} \epsilon^{-rV_\mathcal{F}}
$$
For a universal constant $K$ and $0<\epsilon <1$.

\end{theorem}

Theorem \hyperref[vcsubgraph]{\ref*{vcsubgraph}}, Lemma \hyperref[squared]{\ref*{squared}} and \ref{squared1} actually yield two different bounds for $N_1(\epsilon, \mathcal{L}_{\G}, P)$:

\begin{lemma} \label{subgbound10}
\begin{equation}\label{subgbound2}
N_1(\epsilon, \mathcal{L}_{\G}, P) \leq K(2V_{\G})(16e)^{2V_{\G}} \epsilon^{-2V_{\G}+1}
\end{equation}
\begin{proof}
by Lemma \ref{squared} the class $\mathcal{L}_{\G}$ is a VC subgraph class with $V_{\mathcal{L}_{\G}} \leq 2 V_{\mathcal{G}}-1 $. Put this upper bound and $r=1$ in Theorem \ref{vcsubgraph}.
 
\end{proof}
\end{lemma}

\begin{lemma}
\[
N_1(\epsilon, \mathcal{L}_{\G},P) \leq N_2(\epsilon/2, \mathcal{G}, P) \leq K(V_{\mathcal{G}}-1) (16e)^{V_{\G}-1} \left(\dfrac{\epsilon}{2}\right)^{-2V_{\G}}
\]
\begin{proof}
First inequality comes from lemma \ref{squared1}. Second inequality is Theorem \ref{vcsubgraph} with $r=2$.
\end{proof}
\end{lemma}

\end{subsection}
\begin{subsection}{Parametric Classes}\label{Parametric}
Another large class of classes of functions are the so called \emph{parametric classes} where the functions are indexed by a subset of parameters in $\mathbb{R}^d$. The following is a direct adaptation of \hyperref[vaart98]{\cite{vaart98}, Example 19.7}:

\begin{theorem}\label{vandervaart}
Let $\mathcal{F}:=\{f_\theta:\theta \in \Theta \}$ be a collection of measurable functions indexed by a bounded subset $\Theta \subset \mathbb{R}^d$. Denote by $\| \cdot \|$ the euclidean norm in $\mathbb{R}^d$. Suppose that there exists a measurable function $m$ such that
$$
|f_{\theta_1}(x) - f_{\theta_2}(x)| \leq m(x)\|\theta_1-\theta_2\|
$$
for every $\theta_1, \theta_2 \in \Theta$. If $P(|m|^p) < \infty$, then
$$
N_p(\epsilon, \mathcal{F}, P) \leq K \left (\dfrac{||m||_{p,P}  \diam(\Theta)}{\epsilon}\right )^d
$$
where K is a constant which depends only on $\Theta$ and $d$ and $||m||_{p,P}= \left(P(|m|^p)\right )^{\frac{1}{p}}$.
\end{theorem}

For the proof of this theorem we will need to introduce the concept of \emph{bracketing number} and their relation to covering numbers.

\begin{definition}\label{bracketing}
Given two functions $l$ and $u$ with domain $S$, the bracket $[l,u]$ is the set of all functions $f$ with domain $S$ and with $l \leq f \leq u$. Let $P$ be a probability measure on $S$ and $p\geq 1$. Denote $\|f\|_{p,P}=\left(P(|m|^p)\right )^{\frac{1}{p}}$. A bracket of size $\epsilon$ is a bracket $[l,u]$ with $\| u-l \|_{p,P} < \epsilon $. The bracketing number $N_{[p]}(\epsilon, \mathcal{F}, P)$ is the minimum number of brackets of size $\epsilon$ needed to cover $\mathcal{F}$. l and u need not be in $\mathcal{F}$ but are assumed to have finite $L_p(P)$ norms.
\end{definition}

\begin{lemma}{(\cite{vaart96}, Theorem 2.7.11)}\label{par1} Let $\mathcal{F} = \{ f_t : t \in T \}$ be a class of functions indexed by a metric space $(T,d)$. Supppose that for some fixed function $m$ and every $s,t \in T$:
\[
|f_s(x)-f_t(x) | \leq d(s,t) m(x)
\]
Then for any measure $P$ and $p\geq 1$:
\[
N_{[p]}(2\epsilon \| m\|_{p,P}, \mathcal{F}, P) \leq N(\epsilon, T, d)
\]
\begin{proof}
Let $t_1, \ldots, t_p$ be an $\epsilon$-cover for $T$ with the metric $d$. Then the brackets $[f_{t_i}- \epsilon m, f_{t_i}+\epsilon m ]$ have size $2\epsilon \| m\|_{p,P}$. Let $f_t \in \mathcal{F}$ and $t_i$ such that $d(t,t_i) < \epsilon $. Then $|f_t(x)-f_{t_i}(x)| \leq d(t,t_i) m(x) < \epsilon F(x)$ so that $f_t$ is in the bracket $[f_{t_i}- \epsilon m, f_{t_i}+\epsilon m ]$.
\end{proof}
\end{lemma}

\begin{lemma}\label{relbrack}
\[
N_p(\epsilon, \mathcal{F}, P) \leq N_{[p]}(2\epsilon, \mathcal{F}, P)
\]
\begin{proof}
If $f$ is in the bracket of size $2\epsilon$ $[l,u]$ then it is in the ball of radius $\epsilon$ around $(l+u)/2$.

\end{proof}
\end{lemma}

\textbf{Proof of Theorem \ref{vandervaart}}
\begin{proof} By Lemma \ref{par1} we have that $N_{[p]}(2\epsilon \| m\|_{p,P}, \mathcal{F}, P) $ is bounded by the $\epsilon$ covering number of $\Theta$ in the euclidean metric. After a translation, $\Theta$ is contained in the cube $[0,\diam(\Theta)]^d$. This cube can be covered by $ ( \diam(\Theta)/\epsilon)$ cubes of side $\epsilon$. The circumscribed balls have radius $\sqrt{d}\epsilon/2$ and they still cover the cube. The centers of these balls may be any $x \in \mathbb{R}^d$. As the covering numbers in this case are translation invariant we conclude that 
\begin{gather*}
N(\sqrt{d}\epsilon/2, \Theta \subset \mathbb{R}^d, d) \leq \left( \dfrac{\diam(\Theta)}{\epsilon} \right)^d \\
\end{gather*}
Changing $\epsilon$ for $2\epsilon/\sqrt{d}$ and using Remark \ref{remcov} yields that 
\begin{gather*}
N(\epsilon, \Theta, d) \leq \left ( \dfrac{\sqrt{d}\diam(\Theta)}{\epsilon} \right )^d
\end{gather*}
Now by Lemma \ref{relbrack} and \ref{par1} we get 
\begin{gather*}
N_p(\epsilon  \| m\|_{p,P} , \mathcal{F}, P) \leq N_{[p]}(2\epsilon \| m\|_{p,P}, \mathcal{F}, P) \leq N(\epsilon, \Theta, d)
\end{gather*}
putting things together:
\begin{gather*}
N_p(\epsilon, \mathcal{F}, P) \leq \sqrt{d}^d\left ( \dfrac{\| m\|_{p,P} \diam(\Theta)}{\epsilon} \right )^d
\end{gather*}

\end{proof}

\begin{lemma}\label{vandervaart2} Let $\mathcal{G}$ be a class of positive functions on a set $S$ satisfying the hypothesis of theorem \hyperref[vandervaart]{\ref*{vandervaart}} and bounded by 1, then the class $\mathcal{L}_{\G}$ with domain $S\times [0,1]$ satisfies
$$
N_1(\epsilon, \mathcal{L}_{\G},P) \leq 2^d \sqrt{d}^d \left (\dfrac{||m||_{1,P}  \diam(\Theta)}{\epsilon}\right )^d
$$
\begin{proof} $ \mathcal{L}_{\G} $ is also a parametric class with parameter set $\Theta \subset \mathbb{R}^d$.
\begin{align*}
|(g_{\theta_1}(x)-y)^2 -(g_{\theta_2}(x)-y)^2| &= |(g_{\theta_1}(x)-y)-(g_{\theta_2}(x)-y)||(g_{\theta_1}(x)-y)+(g_{\theta_2}(x)-y)| \\ &\leq 2 |g_{\theta_1}(x)- g_{\theta_2}(x)| \leq 2m(x) \|\theta_1-\theta_2\|
\end{align*}
Plug in $2m(x)$ in \hyperref[vandervaart]{\ref*{vandervaart}} to obtain the result.
\end{proof}
\end{lemma}
As in the section of VC subgraph classes we have another bound arising from Lemma \hyperref[squared1]{\ref*{squared1}}

\begin{lemma}
\[
N_1(\epsilon, \mathcal{L}_{\G}, P) \leq N_2(\epsilon/2, \G, P) \leq 2^d \sqrt{d}^d  \left (\dfrac{||m||_{2,P}  \diam(\Theta)}{\epsilon}\right )^d
\]
\begin{proof} First inequality comes from Lemma \ref{squared1}. The second inequality comes from the proof of Theorem \ref{vandervaart}.

\end{proof}
\end{lemma}

\end{subsection}

\begin{subsection}{Revisiting Vapnik's $\eta$ trick}\label{Trick2}

In section \hyperref[Trick]{\ref*{Trick}} we followed a method to produce bounds for the real risk of a classifier in terms of its empirical risk and a complexity term. This method can also be used to bound the real risk in the context of regression estimation, where $\mathcal{F}$ is a class of positive functions bounded by 1 such that we can bound, for any probability measure, the $L_1(P_n)$ $\epsilon$-covering numbers by a polynomial in $\epsilon^{-1}$. Precisely, we will prove the following result:

\begin{theorem}\label{interval1} Let $\G$ be a class of functions with domain $S$ and range $[0,1]$. Let $(X,Y)$ be a random pair taking values in $S\times[0,1]$ and $(X_1,Y_1), \ldots, (X_n,Y_n)$ be an i.i.d random sample from the distribution of $(X,Y)$. Moreover suppose that 
\[
N_1(\epsilon, \mathcal{F}, P_n) \leq A \epsilon^{-W}
\] 
For some constants $A,W$ depending only on the class $\G$. Let Let $B=8^{W+1}A$, $R_n=n/128$, $Z =W/2$. For $ n \geq 384 Z (B^{-1}\eta)^{\frac{1}{Z}}$ With probability at least $1-\eta$, simultaneously for all $g \in \G$
\[
P(f) \leq P_n(f) + \sqrt{ \dfrac{Z\log(R_n/Z)}{R_n}+ \dfrac{\log(B)}{R_n}-\dfrac{\log(\eta)}{R_n} }
\]

\end{theorem}

To prove this bound we will need the following definition:

\begin{definition}{(Lambert's $\W$ function)} The function $f: [0, \infty) \rightarrow \mathbb{R}$, $x \mapsto x e^x$ is injective and its range is $[0,\infty)$. Denote its inverse function by $\W: [0, \infty) \rightarrow [0, \infty)$. This function satisfies the equation
\[
x= \W(x)e^{\W(x)} \, \, \, \, \mbox{ $x \geq 0$}
\]
\end{definition}

The $\W$ function has two important properties: 
\begin{remark}\label{wproperties}
For $x \geq 3$, $\mathcal{W}(x) \leq \log(x)$. $\lim_{x\to \infty} \mathcal{W}(x)/\log(x) =1$.
\end{remark}

\textbf{Proof of Theorem \ref{interval1}. } 
\begin{proof}
By Theorem \ref{bound} we have
\begin{align*}
P \{||P_n-P|| > \epsilon \} &\leq 8 P \left ( N_1 (\epsilon/8, \mathcal{F}, P_n) \right ) \exp(-n\epsilon^2/128) \\ 
&\leq 8 A \left(\dfrac{\epsilon}{8} \right)^{-W}\exp(-n\epsilon^2/128) \\ & = B \left ( \epsilon^2\right)^{-Z} \exp(-R_n \epsilon^2 )
\end{align*}

Let $\gamma = \epsilon^2$ and equate the previous bound to $\eta$
$$
\eta= B \gamma^{-Z} \exp (-R_n \gamma)
$$

Work out the previous expression for $\eta$:

\begin{align*}
& \eta^{-1}B= \gamma^{Z} \exp (R_n \gamma) \Rightarrow (\eta^{-1}B)^{\frac{1}{Z}}= \gamma \exp \left( \dfrac{R_n \gamma}{Z} \right) \\
& \Rightarrow \dfrac{R_n (\eta^{-1}B)^{\frac{1}{Z}}}{Z} = \dfrac{R_n \gamma}{Z}  \exp \left( \dfrac{R_n \gamma}{Z} \right) \Rightarrow \mathcal{W} \left ( \dfrac{R_n}{Z} (\eta^{-1}B)^{1/Z}  \right)= \dfrac{R_n \gamma}{Z} \\
& \Rightarrow \epsilon = \sqrt{ \dfrac{Z}{R_n}  \mathcal{W} \left ( \dfrac{R_n (\eta^{-1}B)^{1/Z}}{Z}    \right)  }
\end{align*}

By Remark \ref{wproperties}, if ${R_n (\eta^{-1}B)^{1/Z}}    \geq 3Z$ we have that
$$
\epsilon \leq \sqrt{ \dfrac{Z\log(R_n/Z)}{R_n}+ \dfrac{\log(B)}{R_n}-\dfrac{\log(\eta)}{R_n} }
$$
We conclude that with probability at least $1-\eta$ and $n$ large enough so that the previous condition holds, for any $g \in \G$
\[
P(g)-P_n(g) \leq |P(g)-P_n(g)| \leq ||P-P_n|| \leq \epsilon 
\]
\end{proof}
\end{subsection}
\end{section}

\begin{section}{Complexity Penalty Based on Local Metric Entropy}\label{Penalty}
\begin{subsection}{Definitions}\label{definitions}

\begin{definition}Let $\mathcal{G}$ be a class of real valued functions with domain $S$ and range $[0,1]$, and let $P$ be a probability measure on $S$. Let $r\geq 1$ and $B_P(g_0,\delta)$ be a ball of radius $\delta$ around some $g_0 \in \mathcal{G}$ with respect to the $L_r(P)$ metric. 
\[
J(r, \delta, \mathcal{G},g_0, P)= \int_0^\delta \sqrt{\log N_r(u, B_P(g_0,\delta),P)}du
\]
is the local r-metric entropy integral of $\mathcal{G}$. Let $X_1, \ldots, X_n$ be an i.i.d. random sample from the distribution $P$
\[
J(r, \delta, \mathcal{G}, g_0, P_n)=\int_0^\delta \sqrt{\log N_r(u, B_{P_n}(g_0, \delta),P_n)}du
\]
is the random local r-metric entropy integral of $\mathcal{G}$ with respect to the random sample $X_1, \ldots, X_n$.
\[
J(r, \delta, \mathcal{G}, g_0)=\int_0^\delta \sup_{Q}  \sqrt{\log N_r(u, B_Q(g_0,\delta),P)}du
\]
is the uniform local metric entropy integral of $\mathcal{G}$, where the supremum is taken over all probability distributions ons $S$.
\end{definition}
These definitions will be used to define a complexity penalty for the class $\G$.  Recall the following theorem that can be found in \hyperref[pollard84]{\cite{pollard84}, Chapter 2, Theorem 37}.

\begin{theorem}\label{pollardrate}
For each n, let $\mathcal{F}_n$ be a permissible class of positive functions bounded by 1, whose covering numbers satisfy
\begin{equation}\label{condition1}
\sup_Q N_1(\epsilon, \mathcal{F}_n,Q) \leq A \epsilon^{-W} \, \, \, \, \, \mbox{for } 0<\epsilon <1
\end{equation}

with constants $A,W$ not depending on $n$. Let $\alpha_n$ be a non-increasing sequence of positive numbers for which $n\delta_n^2 \alpha_n^2 \gg \log(n)$. If $|f| \leq 1$ and $(Pf^2)^{1/2} \leq \delta_n$ for each $f\in \mathcal{F}_n$, then
\[
\sup_{\mathcal{F}_n} |P_n(f)-P(f)| \ll \delta_n^2 \alpha_n \, \, \, \, \, a.s.
\]
where the notation $a_n \gg b_n $ means $ b_n/a_n \rightarrow 0 $ as $n \to \infty$.

\end{theorem}

For all n, let $ \mathcal{F}_n=\mathcal{L}_{\G}$ for a class of positive functions $\G$ with envelope 1 and $Y$ a $[0,1]$ valued random variable. Suppose that for $\mathcal{L}_{\G}$ the condition \hyperref[condition1]{(\ref*{condition1})} holds. By lemma \hyperref[env1]{\ref*{env1}} the class $\mathcal{L}_{\G}$ also comprise positive functions bounded by 1. For all $n$, let $\delta_n=1$ and $\alpha_n=\log(n)/\sqrt{n}$ then $n \delta_n^2 \alpha_n^2= \log^2 (n) \gg \log(n)$ so the hypothesis of the theorem are satisfied and we get

\[
\sup_{\G} |L_n(g)-L(g)| \ll \log(n)/\sqrt{n} \, \, \, \, \, a.s.
\]

Then by lemma \hyperref[Lbound]{\ref*{Lbound}} $L(\hat{g}_n)- L(g_0) \ll \log(n)/\sqrt{n}$  almost surely, where $g_0$ is the best estimator of the class. Then almost surely the function $(\hat{g}_n-y)^2$ lies within a ball of radius $\delta_n$ in the $L_1$ metric with center $(g_0-y)^2$. \\

For a class $\G$ of positive functions bounded by $1$ defined on a set $S$, an i.i.d. random sample $(X_1,Y_1), \ldots, (X_n,Y_n)$ from the distribution of a random pair $(X,Y)$ taking values in $S\times [0,1]$, $g_0 \in \G$ a function such that $L(g) = \inf_{g \in \G} L(g)$ and $\delta_n= \log(n)/\sqrt{n}$ we will further analyze the quantity
\[
J(2, \delta_n, \mathcal{G}, g_0, P_n) = \int_0^{\delta_n} \sqrt{ \log N_2 (u, B_{P_n} (g_0, \delta_n),P_n))}du
\]

\end{subsection}

\begin{subsection}{Estimating the Local Metric Entropy of a class of functions}\label{estimating}
Let $X_1, \ldots, X_n \in \mathbb{R}^k$ be an i.i.d. sample drawn from an unknown distribution. Suppose the distribution has a density in $\mathbb{R}^k$ with the Lebesgue Measure. Let $\mathcal{G}$ be a class of functions with domain $\mathbb{R}^k$. Define the seminorm $L_2(P_n)$:
\[
\|g\|_n^2 :=\dfrac{1}{n} \sum_{i=1}^n |g(X_i)|^2
\]
And the $L_2(P)$ norm
\[
\|g\|^2 :=P(g^2)
\]

Both induce a pseudometric on $\mathcal{G}$. Suppose $\mathcal{G}:=\left \{ g_{\theta}: \theta \in \Theta \subset \mathbb{R}^d \right \} $ is a family of parametric functions and let $g^*$ be some function in $\mathcal{G}$. Let $B_n^*(\delta)$ be the ball around $g^*$ with radius $\delta$ in the $L_2(P_n)$ metric and let $B^*(\delta)$ be the ball around $g^*$ with radius $\delta$ in the $L_2(P)$ norm. We want to estimate $\log N_2(\epsilon, B_n^*(\delta), P_n)$ and $\log N_2(\epsilon, B^*(\delta), P)$, the  random local metric entropy and the local metric entropy of this class of functions. Throughout we suppose that the map $\theta \mapsto g_\theta$ is bijective.  \\

{\bf Random Local Metric Entropy for Linear Regression }\\

First consider the case where $g_{\theta}(x)= \theta ^Tx$, $\theta \in \mathbb{R}^k$ and $g^*=g_0=0$, then $B_n^*(\delta) = \left \{g_\theta \in \G: \|g_\theta\|_n^2 < \delta^2 \right \} \\$ and define 

\[
\mathbb{X}_n= \dfrac{1}{n} \sum_{i=1}^n X_i X_i^T
\]

Then we have the following equality

\[
\|g_\theta \|_n^2 = \dfrac{1}{n} \sum_{i=1}^n \left |\theta^TX_i \right|^2= \dfrac{1}{n} \sum_{i=1}^n \left ( \theta^TX_i\right ) \left ( \theta^TX_i\right )^T = \dfrac{1}{n} \sum_{i=1}^n  \theta^TX_i X_i^T \theta= \theta^T \mathbb{X}_n \theta 
\]
\\
\begin{lemma} For $n\geq k$,  with probability 1 the matrix $\mathbb{X}_n$ is positive definite and symmetric.
\begin{proof}
For each $i$, the matrix $X_iX_i^T$ is clearly symmetric and it is positive semidefinite as for any nonzero vector $z$, 

\[z^TX_iX_i^Tz=(X_i^Tz)^T(X_i^Tz)=|X_i^Tz|^2 \geq 0 \]
This implies that $\mathbb{X}$ is also symmetric and positive semidefinite. Notice that the equation $|X_i^Tz|=0$ Defines a random hyperplane in $\mathbb{R}^k$, and since the $X_i$'s are independent, the probability that two of such hyperplanes are linearly dependent is zero, implying that for $n \geq k $ the intersection of all such hyperplanes has zero dimension. This is equivalent to $\mathbb{X}_n$ being positive definite with probability one for $n\geq k$.

\end{proof}
\end{lemma}

Let $\mathbb{X}_n=Q^TDQ$ be an orthogonal diagonalization of $\mathbb{X}_n$, $\lambda^T=(\lambda_1, \ldots, \lambda_k) > 0$ its vector of eigenvalues in decreasing order and $(Q\theta)^T=\alpha^T = (\alpha_1, \ldots, \alpha_k)$. Then 

\[
\|g_\theta \|_n^2 = \theta^T Q^T D Q \theta = (Q\theta)^T D (Q\theta) = \sum_{i=1}^n \lambda_i \alpha_i^2= \lambda^T \alpha^2
\]

\[
\dfrac{\|g_\theta \|_n^2}{\delta^2} = \sum_{i=1}^n \lambda_i \dfrac{\alpha_i^2}{\delta^2} = \sum  \dfrac{\alpha_i^2}{\left (\delta/ \sqrt{\lambda_i} \right )^2}
\]

Then the set $\{ \theta : \|g_\theta \|_n^2 < \delta^2 \}$ is actually an ellipsoid in euclidean space with elliptical radii equal to $\delta/ \sqrt{\lambda_i} $.\\

Notice that the following holds
\[
\|g_\theta-g_{\theta'}\|^2 < \delta^2 \Leftrightarrow  (\theta - \theta')^T\mathbb{X}_n(\theta-\theta') < \delta^2
\]
And because $\mathbb{X}_n$ is symmetric and positive definite, it defines an inner product and hence, a metric on the space of parameters. The map $f(g_\theta)= \theta$ is a bijection between the set of functions and the set of parameters, and the conditions of lemma \ref{preservemetric} hold with $K=1$. In the same manner, the map $h(\theta)=Q\theta$ is bijective between the space of parameters $\Theta$ with the metric induced by $\mathbb{X}_n$, and the space of parameters  $Q\Theta $ with the metric induced by the diagonal matrix $D$. \\

Finally the map $h(\alpha)=\sqrt{D}\alpha$ is a bijection between the space $Q \Theta $ and its image, and $\alpha^T D \alpha < \delta^2 \Leftrightarrow h(\alpha)^Th(\alpha) < \delta^2$ so the conditions of lemma \ref{preservemetric} hold with $K=1$ and with the metric of the image space being the usual euclidean metric. Under these maps, the image of the ball  $B_n^*(\delta) = \left \{g_\theta \in \G: \|g_\theta\|_n^2 < \delta^2 \right \} $ is contained in the ball  $ \left \{ \beta \in \mathbb{R}^k: \beta^T\beta < \delta^2 \right \} $ which is a sphere with radius $\delta$ in euclidean space. We have the following results:\\

\begin{lemma} Let $B_\delta^k$ be a ball with radius $\delta$ in euclidean space $\mathbb{R}^k$. In the case of linear regression $g_\theta = \theta^Tx$ and $g^*=\theta^*x$ we have
\[
N_2(u, B_n^*(\delta), P_n) \leq N_2(u, B^k_\delta, d)
\]

\end{lemma}

\begin{lemma}\label{rogers} {\bf (Rogers, \cite{rogers63})}  Let $\mathcal{N}(B_\delta^k)$ be the covering number of $B_\delta^k$ with balls of radius 1 and $k\geq9$. then

\[
\mathcal{N}(B_\delta^k) \leq \left \{  \begin{array}{lr}
        Ck^{5/2} \delta^k  & : \delta < k \\
       Ck \log (k) \delta^k & : \delta \geq k
     \end{array} \right.
\]
Where $C$ is an absolute constant.
\end{lemma}

Note that $\mathcal{N}(u, B_\delta^k)$, the covering number with balls of radius $u$, is equal to $\mathcal{N}(B_{\delta/u}^k)$. We get the following corollary

\begin{corollary}\label{entropybounds}

\[
\log N_2(u, B_n^*(\delta), P_n) \leq \left \{  \begin{array}{lr}
     \log(Ck^{5/2}) + k \log (\delta/u) & : \delta/u < k    \\
      \log(Ck \log k) + k \log (\delta/u) & : \delta/u  \geq k 
     \end{array} \right.
\]
\end{corollary} 

\begin{remark} To ease the computations, we may use the following well known bound for the u-covering number of a ball with radius $\delta$ in euclidean space $\mathbb{R}^k$
\[
N(u, B_\delta^k,d) \leq \left( \dfrac{3\delta}{u} \right )^k
\]
\end{remark}

We conclude the following
\begin{corollary}
\[
N_2(u, B_n^*(\delta), P_n) \leq \left( \dfrac{3\delta}{u} \right )^k
\]
\end{corollary}

{\bf Random Local Metric Entropy for General Linear Regression}\\

Now consider the case where $\mathcal{G}=\{g_\theta = \theta_1 \psi_1(x) + \ldots + \theta_d \psi_d(x) : \theta \in \Theta \subset \mathbb{R}^d\}$, $x \in \mathbb{R}^k$ and $\{ \psi_i \}$ linearly independent functions. We want to estimate again the metric entropy of a ball with respect to the semi norm given by an i.i.d. random sample $X_1, \ldots, X_n$, where the distribution has a density in $\mathbb{R}^k$. In this case
\[
\| g_\theta \|^2 = \dfrac{1}{n} \sum (\theta^T\psi(X_i))^2 =  \dfrac{1}{n} \sum (\theta^T\psi(X_i))(\psi(X_i)^T \theta) = \theta^T \left(\dfrac{1}{n} \sum \psi(X_i)\psi(X_i)^T \right ) \theta 
\]
where $\psi(x)=(\psi_1(x), \ldots, \psi_d(x))^T$. Define $\psi \mathbb{X}_n  = \left(\dfrac{1}{n} \sum \psi(X_i)\psi(X_i)^T \right )$. A similar argument to that of the previous section shows that this matrix is symmetric and for $n \geq d$ it is positive definite with probability one. A straightforward generalization of the previous section yields
\begin{lemma}\label{boundlocal}  For the case of general linear regression 
\[
N_2(u, B_n^*(\delta), P_n) \leq  \left( \dfrac{3\delta}{u} \right )^d
\]

\end{lemma}

{\bf Local Metric Entropy for General Linear Regression}\\

Let $P$ be a probability measure on a set $S$ and $\mathcal{G}=\{g_\theta = \theta_1 \psi_1(x) + \ldots + \theta_d \psi_d(x) : \theta \in \Theta \subset \mathbb{R}^d\}$ for $\psi_i: S \rightarrow \mathbb{R}$ linearly independent functions. We estimate $N_2(u, B^*(\delta), P)$. Let $\psi=( \psi_1, \ldots, \psi_d)$. Note that
\[
P\left((\theta^T \psi)^2\right)= P(\theta^T \psi \psi^T \theta) = \theta^T P(\psi \psi^T) \theta > 0
\]
by linearity and independence. Then the matrix $\Psi= P(\psi \psi^T)$ is symmetric and positive definite. The arguments of the previous section generalize in this case to conclude that, for any probability measure
\[
N_2(u, B^*(\delta), P) \leq  \left( \dfrac{3\delta}{u} \right )^d
\]

Now in all these cases, for $\delta_n= \log(n)/\sqrt{n}$
\[
\int_0^{\delta_n} \sqrt{\log N_2(u, B(\delta_n), P) du} \leq \int_0^{\delta_n} \sqrt{ \log \left( \dfrac{3\delta_n}{u} \right )^d} du
\]
Upon the substitution $3 \delta_n /u$ = $1/z$, this bound  becomes
\[
 3 \delta_n \sqrt{d} \int_0^{1/3} \sqrt{ \log \left( \dfrac{1}{z} \right) } dz = 3A \log(n) \sqrt{\dfrac{d}{n}} 
\]
Where $A= \int_0^{1/3} \sqrt{ \log \left( \dfrac{1}{z} \right) } dz$.

\begin{corollary}\label{corloc}
\[
\int_0^{\delta_n} \sqrt{\log N_2(u, B(\delta_n), P) du} \leq  3A \log(n) \sqrt{\dfrac{d}{n}} 
\]
\end{corollary}

\end{subsection}

\end{section}

\begin{section}{Strong Consistency of a Regression Rule}\label{Strong}
Let $\G_1, \G_2, \ldots $ be a sequence of classes of positive functions on a set $S$ bounded by 1. Let $(X,Y)$ be a random pair taking values in $S\times[0,1]$ and $(X_1,Y_1), \ldots, (X_n,Y_n)$ a random sample from the probability distribution $P$ of $(X,Y)$. Let $\hat{g}_{n,j}$ be the least squares estimator of the class $\G_j$ and let $r(n,j)$ be a function depending on the index $j$ of the class and $n$ the sample size. 

\begin{definition}The Structural Risk of the class $\G_j$ for the complexity penalty $r(n,j)$ is the random variable
\[
\widetilde{L}_{n,j}= \widehat{L}_n(\hat{g}_{n,j})+r(n,j)
\]
Let $g^*_n$ be a function among the $\hat{g}_{n,j}$ minimizing such quantity. We refer to this function as the estimator based on Structural Risk Minimization.
\end{definition}

Let $L^*$ be the Bayes Risk. Assume that
\begin{equation}\label{one}
\inf_{j} \inf_{g \in \G_j} L(g) = L^*
\end{equation}
We want to find sufficient conditions for the complexity penalty $r(j,n)$ such that the estimator based on Structural Risk Minimization is strongly consistent, as defined in \hyperref[consistency]{\ref*{consistency}}.\\

Consider the decomposition
\begin{equation}\label{decomp}
L(g^*_n)-L^*= \left ( L(g^*_n) - \inf_{j \geq 1} \widetilde{L}_{n,j} \right ) + \left( \inf_{j \geq 1} \widetilde{L}_{n,j} - L^*   \right )
\end{equation}
It suffices to show both terms converge to zero almost surely. For the first term we have the following which is a direct adaptation of \hyperref[devroye96]{\cite{devroye96} Theorem 18.2}:
\begin{lemma}\label{firstterm}
\begin{align*}
P \left  \{ L(g^*_n) - \inf_{j \geq 1} \widetilde{L}_{n,j} > \epsilon \right \} \leq \sum_{j=1}^\infty P  \left \{ \sup_{g \in \G_j} |L(g)- \widehat{L}_n(g)| > \epsilon + r(n,j) \right \}
\end{align*}
\begin{proof}
\begin{align*}
P \left  \{ L(g^*_n) - \inf_{j \geq 1} \widetilde{L}_{n,j} > \epsilon \right \} & \leq P \left \{ \sup_{j \geq 1 } \left ( L(\hat{g}_{n,j}) - \widetilde{L}_{n,j} \right ) > \epsilon \right \} \\
& = P \left \{ \sup_{j \geq 1 } \left ( L(\hat{g}_{n,j}) - \widehat{L}_n(\hat{g}_{n,j})-r(n,j) \right ) > \epsilon \right \} \\
& \leq P \left \{ \sup_{j \geq 1 } \left | L(\hat{g}_{n,j}) - \widehat{L}_n(\hat{g}_{n,j})\right | > \epsilon + r(n,j) \right \} \\
& \leq \sum_{j=1}^\infty P \left \{  \left | L(\hat{g}_{n,j}) - \widehat{L}_n(\hat{g}_{n,j})\right | > \epsilon + r(n,j) \right \} \\
& \leq \sum_{j=1}^\infty P  \left \{ \sup_{g \in \G_j} |L(g)- \widehat{L}_n(g)| > \epsilon + r(n,j) \right \}
\end{align*}
\end{proof}
\end{lemma}

Focusing on the second term we have the following lemma which is an adaptation of \hyperref[devroye96]{\cite{devroye96} Theorem 18.2}:
\begin{lemma}\label{second term} Suppose that
\[
\inf_{j} \inf_{g \in \G_j} L(g) = L^*
\]
and that $r(n,j)\rightarrow 0$ for $n \rightarrow \infty$ and any index $j$. Moreover suppose that for any index $j$ and any $\epsilon > 0$ 
\[
\sum_{n=1}^\infty P \left \{ \sup_{g \in \G_j} | \widehat{L}_n(g)-L(g)| > \epsilon \right \} < \infty
\]
Then
\[
\inf_{j \geq 1} \widetilde{L}_{n,j} - L^*   \rightarrow 0 \, \, \, \, a.s.
\]
\begin{proof}
By the hypothesis, for any $\epsilon > 0$ there exists an integer $k$ such that 
\[
\inf_{g \in \G_k} L(g)-L^* \leq \epsilon
\]
Fix such $k$. It suffices to show that
\[
\limsup_{n \to  \infty} \inf_{j \geq 1 } \widetilde{L}_{n,j} - \inf_{g \in \G_k} L(g) \leq 0 \, \, \, \, a.s.
\]
By hypothesis there exist $n_0$ large enough such that $r(n,k) \leq \epsilon/2$ for $n \geq n_0$. Then for such n
\begin{align*}
P \left \{ \inf_{j \geq 1} \widehat{L}_{n,j} - \inf_{g \in \G_k} L(g) > \epsilon \right \} & \leq P \left \{ \widehat{L}_n(\hat{g}_{n,k})+r(k,n)- \inf_{g \in \G_k} L(g) > \epsilon \right \}  \\
& \leq P \left \{ \widehat{L}_n(\hat{g}_{n,k})- \inf_{g \in \G_k} L(g) > \epsilon/2 \right \} \\
& \leq P \left \{ \sup_{g \in \G_k} | \widehat{L}_n(g)-L(g) | > \epsilon/2 \right \}
\end{align*}
By hypothesis such bound is summable in $n$. By the Borel-Cantelli Lemma we have the desired result.

\end{proof}

\end{lemma}

\begin{subsection}{Structural Risk Minimization for VC subgraph classes}
With these results we are ready to prove the following

\begin{theorem}\label{result1} 
Let $\G_1, \G_2, \ldots $ be classes of functions such that for all $j$ there exists a quantity $W_j \geq 1$ and absolute constants $K,A$ such that for all $0 < \epsilon <1$
\begin{equation}\label{two}
P  \left \{ \sup_{g \in \G_j} |L(g)- \widehat{L}_n(g)| > \epsilon \right \}   \leq K W_j \left ( \dfrac{A}{\epsilon} \right )^{W_j} \exp{(-n \epsilon^2/128)}
\end{equation}
Define the complexity penalty
\[
r(n,j) = \sqrt{\dfrac{128W_j \log(Aen)}{n}}
\]
and suppose that
\begin{equation}\label{three}
\Delta= \sum_{j=1}^\infty \exp(-W_j/2) < \infty
\end{equation}
Let $(X,Y)$ be a random pair such that
\[
\inf_{j} \inf_{g \in \G_j} L(g) = L^*
\]
Where $L^*$ is the Bayes Risk. Then the estimator $g^*_n$ based on Structural Risk Minimization is strongly consistent.

\begin{proof} 

Define the following:
\begin{gather*}
T_1= W_j \left (\dfrac{A}{ r(n,j)} \right )^{W_j} \\
 T_2= \exp(-nr(n,j)^2/128)
\end{gather*}

Then we have that
\begin{align*}
P \left \{ \sup_{g \in \G_j} |L(g)-\widehat{L}_n(g) | > \epsilon + r(n,j) \right \} 
&\leq K W_j \left (\dfrac{A}{\epsilon + r(n,j)} \right )^{W_j} \exp(-n(\epsilon+r(n,j))^2/128)\\
& \leq KW_j  \left (\dfrac{A}{r(n,j)} \right )^{W_j} \exp(-nr(n,j)^2/128) \exp(-n\epsilon^2/128) \\
& = K T_1 T_2  \exp(-n\epsilon^2/128)
\end{align*}

Working out the terms $T_1,T_2$ we have

\begin{align*}
T_1 &= \exp( \log(W_j) + W_j \log(A) - W_j \log(r(n,j)))\\
T_2 &= \exp ( -W_j\log(A) -W_j\log(n)-W_j) \\
T_1T_2 &= \exp( \log(W_j) - W_j \log(r(n,j)) -W_j\log(n)-W_j) \\
T_1T_2& \leq  \exp(-W_j\log(n r(n,j))) 
\end{align*}
Where the last inequality comes from the fact that $\log(W_j)-W_j \leq 0 $. We can bound the last inequality in the following way
\begin{align*}
\log(n r(n,j))&=\log\left(\sqrt{128W_jn\log(Aen)} \right ) \\
& = \dfrac{1}{2}\log\left(128W_jn\log(Aen) \right ) \\
& \geq \dfrac{1}{2}
\end{align*}
For $n \geq A^{-1}$ which does not depend on $j$. In summary we have
\begin{align*}
P \left \{ \sup_{g \in \G_j} |L(g)-\widehat{L}_n(g) | > \epsilon + r(n,j) \right \} & \leq  K T_1 T_2  \exp(-n\epsilon^2/128) \\
& \leq K \exp(-W_j/2)  \exp(-n\epsilon^2/128)
\end{align*}

By lemma \ref{firstterm} and (\ref{three}) we get for $n \geq A^{-1}$

\begin{align*}
P \left  \{ L(g^*_n) - \inf_{j \geq 1} \widetilde{L}_{n,j} > \epsilon \right \} &\leq \sum_{j=1}^\infty P  \left \{ \sup_{g \in \G_j} |L(g)- \widehat{L}_n(g)| > \epsilon + r(n,j) \right \} \\
& \leq K  \sum_{j=1}^\infty \exp(-W_j/2)  \exp(-n\epsilon^2/128) < \infty \\
& = K  \Delta  \exp(-n\epsilon^2/128)
\end{align*}
By the Borel-Cantelli Lemma we get that 
\[
\left ( L(g^*_n) - \inf_{j \geq 1} \widetilde{L}_{n,j} \right ) \rightarrow 0 \, \, \, \, a.s.
\]
Finally by the definition of the complexity term, for any index $j$, $r(n,j) \rightarrow 0$ as $n \to \infty$ and also by (\ref{two})
\begin{align*}
\sum_{n=1}^\infty P \left \{ \sup_{g \in \G_j} | \widehat{L}_n(g)-L(g)| > \epsilon \right \} \leq \sum_{n=1}^\infty K W_j \left ( \dfrac{A}{\epsilon} \right )^{W_j} \exp{(-n \epsilon^2/128)} < \infty
\end{align*}
By lemma \ref{second term} we conclude that also 
\[
\inf_{j \geq 1} \widetilde{L}_{n,j} - L^*   \rightarrow 0 \, \, \, \, a.s.
\]
Then by the decomposition of equation (\ref{decomp}) we conclude that the estimator is strongly consistent.

\end{proof}

\end{theorem}

\begin{corollary} \label{vcsubfsrm} \textbf{Structural Risk Minimization for VC subgraph classes: } Let $\G_1, \G_2, \ldots $ be VC subgraph classes of functions with domain $S$ and range $[0,1]$. By lemma \ref{subgbound10}
\begin{align*}
N_1(\epsilon, \mathcal{L}_{\G_j}, P) &\leq K(2V_{\G_j})(16e)^{2V_{\G_j}} \epsilon^{-2V_{\G_j}} \\
& = K(2V_{\G_j})\left( \dfrac{16e}{\epsilon} \right )^{2V_{\G_j}} 
\end{align*}
This implies by Theorem \ref{importantbound}
\[
P \left \{ \sup_{g \in \mathcal{G}} |L(g)-\widehat{L}_n(g)| > \epsilon \right \} \leq 8 K(2V_{\G_j})\left ( \dfrac{128e}{\epsilon} \right )^{2V_{\G_j}}
 \exp(-n\epsilon^2/128)
\]
Suppose that 
\begin{equation*}
\Delta= \sum_{j=1}^\infty \exp(-V_{\G_j}) < \infty
\end{equation*}
Define the complexity penalty
\[
r(n,j) = \sqrt{\dfrac{256V_{\G_j} \log(128e^2n)}{n}}
\]
Let $(X,Y)$ be a random pair such that
\[
\inf_{j} \inf_{g \in \G_j} L(g) = L^*
\]
Where $L^*$ is the Bayes Risk. Then the estimator $g^*_n$ based on Structural Risk Minimization is strongly consistent.

\end{corollary}

The conditions imposed on the class of functions are not so restrictive. We could take a nested sequence $\G_1 \subset \G_2 \subset \ldots$ of VC subgraph classes such that the VC subgraph index increases by one for each class. The union of all such classes is a class with infinite VC dimension so that $g^*$, the Bayes regression function, most likely belongs to one of the classes, implying the last condition.
\end{subsection}

\begin{subsection}{Structural Risk Minimization in parametric classes}
By lemma \hyperref[vandervaart2]{\ref*{vandervaart2}} we have that for a parametric class with $d$ parameters satisfying the conditions of theorem \ref{vandervaart}:
\[
N_1(\epsilon, \mathcal{L}_{\G},P) \leq 2^d\sqrt{d}^d \left (\dfrac{||m||_{1,P}  \diam(\Theta)}{\epsilon}\right )^d
\]

\begin{example}\label{parametric}
Consider the case of a sequence of linearly independent positive functions on a set $S$, $\psi_1, \psi_2, \ldots$  bounded by 1. Let $\psi^{(j)}=(\psi_1, \ldots, \psi_j)^T$ and consider the class
\[
\G_j= \{ \theta^T \psi^{(j)} : \theta \in \Theta_j \}
\]
Where $\Theta_j= \{ \theta \in \mathbb{R}^j: \sum \theta_i \leq 1 : 0 \leq \theta_i \}$. One can check all this sets have diameter $\sqrt{2}$.
\end{example}

\begin{lemma} $\diam(\Theta_1)=1$ and for $j \geq 2$, $\diam(\Theta_j) = \sqrt{2}$.
\begin{proof}
The assertion for $\Theta_1$ is trivial. For $j\geq 2$ note that $\Theta_j$ is a closed convex polytope. Its diameter is equal to the maximum distance between any two vertices. The vertices of such polytope are the zero vector and the unitary vectors $\{ e_i : 1 \leq i \leq j \}$.

\end{proof}

\end{lemma}

\begin{lemma}
 Let $m_j(x)=\sqrt{j}$, then
\[
|\theta^T_1 \psi^{(j)}(x) - \theta^T_2 \psi^{(j)}(x) | \leq m_j(x) \|\theta^T_1 - \theta^T_2  \|
\]
The bound is tight if for some $x$, $\psi^{(j)}(x)=(1, \ldots, 1)$.
\begin{proof}
\begin{align*}
|\theta^T_1 \psi^{(j)} - \theta^T_2 \psi^{(j)}(x) |  & = |(\theta^T_1 - \theta^T_2) \psi^{(j)}(x) |  
\\ & \leq \| \theta^T_1 - \theta^T_2 \|  \| \psi^{(j)}(x) \| \leq \sqrt{j} \| \theta^T_1 - \theta^T_2 \|
\end{align*}

\end{proof} 
Now let $x$ be such that $\psi^{(j)}(x)=(1, \ldots, 1)$, $\theta^T_1= (\frac{1}{j}, \ldots, \frac{1}{j})$ and $\theta^T_2= (0, \ldots, 0)$. Then
\begin{align*}
|\theta^T_1 \psi^{(j)}(x) - \theta^T_2 \psi^{(j)}(x) | &=|\theta^T_1 \psi^{(j)}(x)| \\
&= \sum_{i=1}^j  \frac{1}{j}=1
\end{align*}
On the other hand
\[
 \sqrt{j} \| \theta^T_1 - \theta^T_2 \| = \sqrt{j} \| \theta^T_1 \| =  \sqrt{  \sum_{i=1}^j  j\frac{1}{j^2}}=1 
\]
\end{lemma}

We would like to obtain an analogous result to Theorem \ref{result1} in this case. 
\begin{theorem}\label{result2}
Let $\G_1, \G_2, \ldots $ be classes of functions such that for all $j$ there exists a quantity $M_j \geq 1$  such that for all $0 < \epsilon <1$
\begin{equation}\label{two2}
P  \left \{ \sup_{g \in \G_j} |L(g)- \widehat{L}_n(g)| > \epsilon \right \}  \leq K (2^j)\sqrt{j}^j \left (\dfrac{M_j}{\epsilon}\right )^j  \exp{(-n \epsilon^2/128)} 
\end{equation}
Define the complexity penalty
\[
r(n,j) = \sqrt{\dfrac{128 j \log(2j^{1/2}M_j n)}{n}}
\]

Let $(X,Y)$ be a random pair such that
\[
\inf_{j} \inf_{g \in \G_j} L(g) = L^*
\]
Where $L^*$ is the Bayes Risk. Then the estimator $g^*_n$ based on Structural Risk Minimization is strongly consistent.

\begin{proof}
Define the following:
\begin{gather*}
T_1= 2^j\sqrt{j}^j \left (\dfrac{M_j}{r(n,j)}\right )^j\\
 T_2= \exp(-nr(n,j)^2/128)
\end{gather*}

Then we have that
\begin{align*}
P \left \{ \sup_{g \in \G_j} |L(g)-\widehat{L}_n(g) | > \epsilon + r(n,j) \right \} 
&\leq  K (2^j)\sqrt{j}^j \left (\dfrac{M_j}{\epsilon+r(n.j)}\right )^j \exp(-n(\epsilon+r(n,j))^2/128)\\
& \leq K(2^j)\sqrt{j}^j \left (\dfrac{M_j}{r(n,j)}\right )^j \exp(-nr(n,j)^2/128) \exp(-n\epsilon^2/128) \\
& =  K T_1 T_2  \exp(-n\epsilon^2/128)
\end{align*}

Working out the terms $T_1,T_2$ we have
\begin{align*}
T_1 &= \exp( j \log(2) + \frac{1}{2} j \log(j) + j \log(M_j) - j \log(r(n.j))    )\\
T_2 &= \exp ( - j \log(2) - \frac{1}{2} j \log(j) - j \log(M_j) -j \log(n)) \\
T_1T_2 & \leq \exp( -j \log(n r(n,j)) ) \\
\end{align*}

We can bound the last inequality in the following way, for $n,j \geq 1$
\begin{align*}
\log(n r(n,j))&=\log\left(\sqrt{128j n\log(2 j^{1/2} M_j n)} \right )  \geq 1
\end{align*}

In summary we have
\begin{align*}
P \left \{ \sup_{g \in \G_j} |L(g)-\widehat{L}_n(g) | > \epsilon + r(n,j) \right \} & \leq  KT_1 T_2  \exp(-n\epsilon^2/128) \\
& \leq K \exp(-j)  \exp(-n\epsilon^2/128)
\end{align*}
By lemma \ref{firstterm} and we get:

\begin{align*}
P \left  \{ L(g^*_n) - \inf_{j \geq 1} \widetilde{L}_{n,j} > \epsilon \right \} &\leq \sum_{j=1}^\infty P  \left \{ \sup_{g \in \G_j} |L(g)- \widehat{L}_n(g)| > \epsilon + r(n,j) \right \} \\
& \leq K  \sum_{j=1}^\infty \exp(-J)  \exp(-n\epsilon^2/128) < \infty \\
\end{align*}

Because $\Delta=\sum_{j=1}^\infty \exp(-j) < \infty$, by the Borel-Cantelli Lemma we get that 
\[
\left ( L(g^*_n) - \inf_{j \geq 1} \widetilde{L}_{n,j} \right ) \rightarrow 0 \, \, \, \, a.s.
\]
Finally by the definition of the complexity term, for any index $j$, $r(n,j) \rightarrow 0$ as $n \to \infty$ and also by (\ref{two2})
\begin{align*}
\sum_{n=1}^\infty P \left \{ \sup_{g \in \G_j} | \widehat{L}_n(g)-L(g)| > \epsilon \right \} \leq \sum_{n=1}^{\infty} K(2^j)\sqrt{j}^j \left (\dfrac{M_j}{\epsilon}\right )^j  \exp{(-n \epsilon^2/128)} 
 < \infty
\end{align*}
By lemma \ref{second term} we conclude that also 
\[
\inf_{j \geq 1} \widetilde{L}_{n,j} - L^*   \rightarrow 0 \, \, \, \, a.s.
\]
Then by the decomposition of equation (\ref{decomp}) we conclude that the estimator is strongly consistent.
\end{proof}

\end{theorem}

\begin{corollary} In the case of example \ref{parametric}
\begin{align*}
P  \left \{ \sup_{g \in \G_j} |L(g)- \widehat{L}_n(g)| > \epsilon \right \}  & \leq 8 P N_1(\epsilon/8, \mathcal{L}_\mathcal{G}, P_n) \exp(-n\epsilon^2/128)
\\ &\leq 8 (2^j)\sqrt{j}^j \left (\dfrac{8\sqrt{2j}}{\epsilon}\right )^j  \exp{(-n \epsilon^2/128)} 
\end{align*}

Define the complexity penalty
\[
r(n,j) = \sqrt{\dfrac{128j \log(16 \sqrt{2}j n)}{n}}
\]

Let $(X,Y)$ be a random pair such that
\[
\inf_{j} \inf_{g \in \G_j} L(g) = L^*
\]
Where $L^*$ is the Bayes Risk. Then the estimator $g^*_n$ based on Structural Risk Minimization is strongly consistent.
\end{corollary}
\end{subsection}

\textbf{Further considerations: }\\

It is not easy to obtain bounds for the local metric entropy of a class of functions $\G$, we have done so in the case of general linear regression which is an important but simple case. In further work we aim to establish similar kind of bounds for example in the case of VC subgraph classes. For a sequence of classes $\G_1, \G_2, \ldots, $, and a probability measure $P$ we aim to find a bound
\[
N_2(u, B(g_{0,j}, \delta), P) \leq \Psi(u, j, \delta)
\]
where $B(g_{0,j}, \delta) \subset \G_j$ is the ball around some $g_{0,j} \in \G_j$ with respect to the $L_2(P)$ metric. And study the complexity penalty 
\[
r(n,j)=\int_0^{\delta_n} \sqrt{\log(\Psi(u,j,\delta_n))}du
\]
for $\delta_n= \log(n)/n$. In the case of \ref{parametric}, By lemma \ref{corloc}, we would define the complexity penalty
\[
\int_0^{\delta_n} \sqrt{\log N_2(u, B(g_{0,j}, \delta_n), P) }du \leq  3A \log(n) \sqrt{\dfrac{j}{n}} =r(j,n)
\]
We would get the following bound
\begin{align*}
P \left \{ \sup_{g \in \G_j} |L(g)-\widehat{L}_n(g) | > \epsilon + r(n,j) \right \} 
&\leq  8 (2^j)\sqrt{j}^j \left (\dfrac{ 8  \sqrt{2j}}{\epsilon+r(n.j)}\right )^j \exp(-n(\epsilon+r(n,j))^2/128)\\
& \leq 8(2^j)\sqrt{j}^j \left (\dfrac{ 8  \sqrt{2j}}{r(n,j)}\right )^j \exp(-nr(n,j)^2/128) \exp(-n\epsilon^2/128) \\
& =  8 T_1 T_2  \exp(-n\epsilon^2/128)
\end{align*}
\end{section}
where 
\begin{align*}
T_1&= 2^j\sqrt{j}^j \left (\dfrac{ 8  \sqrt{2j}}{r(n,j)}\right )^j\\
 T_2&= \exp(-nr(n,j)^2/128)
 \end{align*}
 Working out the terms $T_1,T_2$ we have
\begin{align*}
T_1 &= \exp( j \log(2) + \frac{1}{2} j \log(j) + j \log(8  \sqrt{2j}) - j \log(r(n.j))    )\\
T_2 &= \exp ( - j 9A^2 \log^2(n) /128)) \\
\end{align*}
Now $ \log(r(n.j)) \leq \log(2A \sqrt{j})$ so that $-j \log(r(n,j)) \geq - j \log(3A \sqrt{j})$. Then
\begin{align*}
T_1T_2 \geq \exp \left( j  \left [ \log\left( \dfrac{2 \sqrt{j} 8 \sqrt{2j} }{ 3A \sqrt{j}  } \right) - \dfrac{9A^2 \log^2(n)} {128} \right ] \right) = \exp \left( j  \left [ \log\left( \dfrac{16 \sqrt{2j} }{ 3A   } \right) - \dfrac{9A^2 \log^2(n)} {128} \right ] \right) \geq \exp(j)
\end{align*}
For $j$ large enough. Then the obtained bound for 
\[
\sum_{j=1}^\infty P  \left \{ \sup_{g \in \G_j} |L(g)- \widehat{L}_n(g)| > \epsilon + r(n,j) \right \} 
\]
fails to converge. We are not able to mimic the proof of \ref{result2}.

\end{document}